\documentclass[11pt, a4paper, twoside, reqno]{amsart}

\usepackage{charter}
\usepackage{dutchcal}
\usepackage{ragged2e}
\usepackage{enumitem, amsmath, amsthm, amssymb, amsfonts, mathtools, fullpage, ulem, graphicx}
\usepackage{amsbsy, amssymb, amscd, amsfonts, latexsym, amstext, delarray, amsmath, graphicx, color, caption, blkarray, mathtools}
\usepackage{tikz}
\usetikzlibrary{cd}
\usetikzlibrary{matrix, arrows}
\usepackage{graphicx}
\usepackage{euscript}
\usepackage{thm-restate}
\usepackage[colorlinks=true, citecolor=black, linkcolor=black, urlcolor=blue]{hyperref}
\newcommand{\seqnum}[1]{\href{http://oeis.org/#1}{\underline{#1}}}
\newcommand{\lucky}{\mathrm{lucky}}
\newcommand{\unlucky}{\mathrm{unlucky}}
\newcommand{\bump}{\mathrm{bump}}

\usepackage{cleveref}
\newtheorem{theorem}{Theorem}[section]
\newtheorem{lemma}[theorem]{Lemma}
\newtheorem{corollary}[theorem]{Corollary}
\newtheorem{proposition}[theorem]{Proposition}
\newtheorem{question}[theorem]{Question}
\theoremstyle{definition}
\newtheorem{definition}[theorem]{Definition}
\newtheorem{example}[theorem]{Example}
\newtheorem{claim}{Claim}
\newtheorem{remark}{Remark}

\newcommand{\Out}{\mathcal{O}}
\newcommand{\Sym}{\mathfrak{S}}
\newcommand{\PF}{\mathrm{PF}}
\newcommand{\MPF}{\mathrm{MPF}}
\newcommand{\BPF}{\mathrm{BPF}}
\newcommand{\NN}{\mathbb{N}}
\newcommand{\N}{\texttt{N}}
\newcommand{\E}{\texttt{E}}
\newcommand{\D}{\texttt{D}}
\newcommand{\U}{\texttt{U}}
\newcommand{\tH}{\texttt{H}}
\newcommand{\calD}{\mathcal{D}}

\newcommand{\calS}{\mathcal{S}}

\title{{ PARKING FUNCTIONS INSPIRED BY BOARD GAMES }}
\date{\today}

\author{Alexander
Clifton}
\address[A.~Clifton]{Department of Theoretical Computer Science, Czech Technical University in Prague, Prague, Czechia}
\email{\textcolor{blue}{\href{mailto:alexander.clifton@fit.cvut.cz}{alexander.clifton@fit.cvut.cz}}}

\author{Pamela
E.
Harris}
\address[P.~E.
Harris]{Department of Mathematical Sciences, University of Wisconsin-Milwaukee, Milwaukee, WI 53211 USA}
\email{\textcolor{blue}{\href{mailto:peharris@uwm.edu}{peharris@uwm.edu}}}
\author{Maryam
Khaqan}
\address[M.~Khaqan]{Department of Mathematics, University of Toronto, Toronto, ON, Canada M5T 3J1}
\email{\textcolor{blue}{\href{mailto:maryam.khaqan@utoronto.ca}{maryam.khaqan@utoronto.ca}}}

\date{}

\begin{document}
	\begin{abstract}
		We introduce two new parking schemes inspired by the classic board games Monopoly and backgammon. We consider a one-way street with $n$ parking spots and a queue of $n$ cars waiting to park. Each car has a preferred parking spot on the street and enters in sequence to park in that spot if it is available. Under the Monopoly parking scheme, if a car finds its preferred spot occupied, it returns to the start of the street and parks in the first available spot it encounters. Under the backgammon parking scheme, if a car finds its preferred spot occupied, it bumps the already parked car out of that spot and takes its place. The bumped car then returns to the start of the street and parks in the first available spot it encounters.

		A preference list is a parking function for a given parking scheme if all cars successfully park under that scheme. We prove that all preference lists are parking functions under both the Monopoly and backgammon parking schemes. However, the final order in which the cars park, known as the outcome of a preference list under a parking scheme, is a permutation of $[n]$ and differs significantly under the Monopoly and backgammon parking schemes.

		For a fixed permutation $\pi$, we provide a product formula for the number of Monopoly parking functions that result in the outcome $\pi$. Analogously, we show that the number of backgammon parking functions that result in the outcome $\pi$ is bounded below by $2^{n-1}$ and bounded above by the $n$-th Catalan number. Furthermore, for several families of outcomes, we provide exact formulas involving Motzkin numbers for the corresponding number of backgammon parking functions.

		Under the Monopoly and backgammon parking schemes, we define the displacement of a car as the {absolute difference} between the spot where that car ultimately parks and its preferred spot. We study displacement multisets that encode the displacements for cars with nonzero displacement. We provide formulas for the number of backgammon parking functions with a fixed displacement multiset of cardinality at most $2$, as well as for the number of Monopoly parking functions whose displacement multiset is a given singleton set. We conclude with some directions for further research.
	\end{abstract}

	\begingroup \def\uppercasenonmath#1{} %
	\let\MakeUppercase\relax %
	\maketitle

	\section{Introduction}
	Consider a one-way street with $n\in \NN=\{1,2,3,\ldots\}$ sequentially numbered parking spots, and a queue of $n$ cars waiting to park in those spots. Each car has a preferred parking spot, and we collect the preferences in a list $\alpha =(a_{1},a_{2},\ldots,a_{n})\in[n]^{n}$, where $a_{i}$ denotes the the preferred parking spot of car $i$ for each $i\in[n] =\{1,2,\dots,n\}$. We call the list $\alpha$ a \textit{preference list}. The cars proceed to park under the \textit{classical parking scheme}: Cars enter sequentially, from car $1$ through car $n$, in that order. Car~$i$ attempts to park in its preferred parking spot $a_{i}$. If the spot is unoccupied, car $i$ parks there, and then car $i+1$ enters the street, attempting to park. If car $i$ finds its preferred parking spot $a_{i}$
	occupied, then car $i$ continues down the one-way street and attempts to park in the first available spot it encounters past its preference. If no such spot exists, meaning all of the spots from $a_{i}$ to $n$ are occupied, then car $i$ exits the street, and parking fails. If under this parking scheme all cars are able to park, then the preference list $\alpha$ is called a \textit{(classical) parking function of length} $n$. For example, $(1,4,4,3,2,2)$ is a parking function of length $6$, resulting in the cars parking as illustrated in \Cref{fig:example of pf}. In contrast, $(2,4,4,3,2,5)$ is not a parking function since no car prefers parking spot $1$ and thus at least one car fails to park.

	Given a classical parking function, the order the cars park on the street is a permutation of $[n]$. We let $\Sym_{n}$ denote the set of permutations of $[n]$, which we write in one-line notation as $\pi=\pi_{1}\pi_{2}\cdots\pi_{n}$. Then, for a parking function $\alpha=(a_{1},a_{2},\ldots,a_{n})\in[n]^{n}$, we let $\Out(\alpha)\in\Sym_{n}$ be the permutation corresponding to the \textit{outcome of $\alpha$}; that is, we have $\Out(\alpha)=\pi=\pi_{1}\pi_{2}\cdots\pi_{n}\in\Sym_{n}$ where $\pi_{j}=k$ means that car $k$ parks in parking spot $j$. Note that the outcome of the example in \Cref{fig:example of pf} is the permutation $\pi=154236\in\Sym_{6}$.
	\begin{figure}
		[h!] \centering \resizebox{\textwidth}{!}{ \begin{tikzpicture}\node at(0,0){\includegraphics[width=1in]{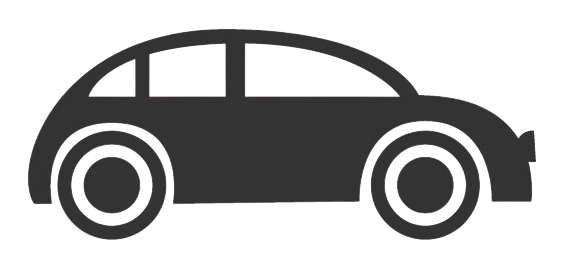}}; \node at(3,0){\includegraphics[width=1in]{car.png}}; \node at(6,0){\includegraphics[width=1in]{car.png}}; \node at(9,0){\includegraphics[width=1in]{car.png}}; \node at(12,0){\includegraphics[width=1in]{car.png}}; \node at(15,0){\includegraphics[width=1in]{car.png}}; \node at(0,-.05){\textcolor{white}{\textbf{1}}}; \node at(3,-.05){\textcolor{white}{\textbf{5}}}; \node at(6,-.05){\textcolor{white}{\textbf{4}}}; \node at(9,-.05){\textcolor{white}{\textbf{2}}}; \node at(12,-.05){\textcolor{white}{\textbf{3}}}; \node at(15,-.05){\textcolor{white}{\textbf{6}}}; \node at(.9,1.){\includegraphics[width=.75in]{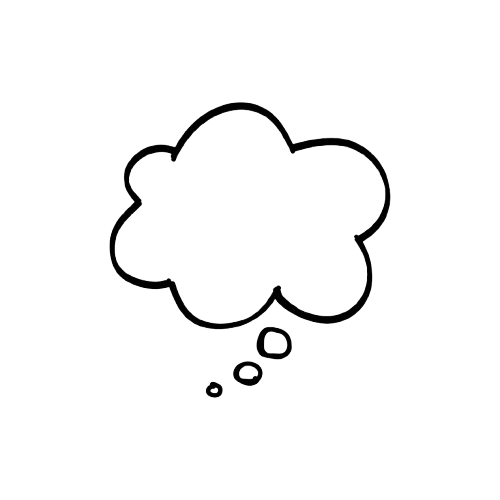}}; \node at(3.9,1.){\includegraphics[width=.75in]{callout.png}}; \node at(6.9,1.){\includegraphics[width=.75in]{callout.png}}; \node at(9.9,1.){\includegraphics[width=.75in]{callout.png}}; \node at(12.9,1.){\includegraphics[width=.75in]{callout.png}}; \node at(15.9,1.){\includegraphics[width=.75in]{callout.png}}; \node at(.9,1.1){\textbf{1}}; \node at(3.9,1.1){\textbf{2}}; \node at(6.9,1.1){\textbf{3}}; \node at(9.9,1.1){\textbf{4}}; \node at(12.9,1.1){\textbf{4}}; \node at(15.9,1.1){\textbf{2}}; \draw[ultra thick](-1.25,-.6)--(1.25,-.6); \draw[ultra thick](1.75,-.6)--(4.25,-.6); \draw[ultra thick](4.75,-.6)--(7.25,-.6); \draw[ultra thick](7.75,-.6)--(10.25,-.6); \draw[ultra thick](10.75,-.6)--(13.25,-.6); \draw[ultra thick](13.75,-.6)--(16.25,-.6); \node at(0,-1){\textbf{1}}; \node at(3,-1){\textbf{2}}; \node at(6,-1){\textbf{3}}; \node at(9,-1){\textbf{4}}; \node at(12,-1){\textbf{5}}; \node at(15,-1){\textbf{6}};\end{tikzpicture} }
		\caption{The
		parking
		outcome
		for
		the
		preference
		list
		$(1,4,4,3,2,2)$
		under
		the
		classical
		parking
		scheme.
		Hence,
		$\Out((1,4,4,3,2,2))=154236$.}
		\label{fig:example of pf}
	\end{figure}

	Parking functions were introduced and studied by Konheim and Weiss in their study of hash functions~\cite{Konheim1966AnOD}. They showed that, if $\PF_{n}$ denotes the set of (classical) parking functions of length $n$, then $|\PF_{n}|=(n+1)^{n-1}$. Since their introduction, parking functions have received considerable attention in the literature due to their numerous connections with other areas of mathematics. For example, the study of parking functions is related to hyperplane arrangements, counting facets of the permutohedron, the QuickSort algorithm, Boolean intervals in the weak Bruhat order of the symmetric group, and the Tower of Hanoi game~\cite{Aguillon2022OnPF,bradt2024unit,unit_perm,elder2024parking,Harris2023LuckyCA,Stanley1996HyperplaneAI}. There are many generalizations of
	classical parking functions often involving either a change to the cars' lengths, their preferences, or the parking scheme; see for examples~\cite{bib:ParkingCompletions,bib:NaplesPF,bib:Colaric2020IntervalPF,bib:ColmenarejoHarris,bib:EhrenborgHapp, bib:SamSpiro}. For more information on parking functions, we recommend the chapter by Yan~\cite{yan2015parking} and for a list of open problems in this area, we refer the interested reader to~\cite{bib:PFCYOA}.

	In studying parking functions, a question of interest is the following: \textit{Supposing that one observes only the final parking outcome and not the preference list, how many distinct preference lists could have resulted in that parking outcome?} In the classical parking function scenario, the result~\cite[Proposition 3.1]{bib:ColmenarejoHarris} is a product of lengths of certain longest contiguous subsequences in the permutation describing the parking outcome. This question has also been considered under other parking schemes including MVP parking functions, $k$-cascading parking functions, $k$-Naples parking functions, pullback parking functions, and $\ell$-interval parking functions~\cite{aguilarfraga2024interval,adenbaum2026parkingcascadessimplestsequence,bib:ColmenarejoHarris,bib:pullbackpfs,MVP}.

	A second area of interest is the study of parking functions with a given \textit{displacement multiset}, formally defined in \Cref{def:displacementvector}. In the classical setting, the displacement of a car $i$ is defined as the difference between the spot car $i$ parks in and the spot it originally preferred. In this paper, we will introduce new parking schemes in which a car may end up parking in a spot earlier than its preferred one, making the aforementioned difference negative. To address this, we henceforth define the displacement $d_{i}$ of car $i$ to be the \textbf{absolute value} of the difference between the spot car $i$ prefers and where it ultimately parks. Note that this definition of displacement agrees with that of displacement under the classical parking scheme, since
	in that setting the difference between a car's final parking spot and its preference is always nonnegative.

	We define the \textit{displacement multiset} to be the multiset consisting of the nonzero values $d_{i}$, and we denote this by $\calD(\alpha)$. Then, given a (classical) parking function $\alpha$, the displacement of $\alpha$, denoted $D(\alpha)$ is the sum of the displacements among all cars. Namely, $D(\alpha)=\sum_{i=1}^{n}d_{i}$, which is equal to the sum of the entries in $\calD(\alpha)$.

	A feature of classical parking functions is that a rearrangement of a parking function is again a parking function. This follows from the inequality characterization of classical parking functions: a preference list $\alpha=(a_{1},a_{2},\ldots,a_{n})\in[n]^{n}$ is a classical parking function if and only if its weakly increasing rearrangement denoted $\alpha^{\uparrow}=(a_{1}',a_{2}',\ldots,a_{n}')$, obtained by sorting the entries into weakly increasing order, satisfies $a_{i}'\leq i$ for all $i\in[n]$. Moreover, given a classical parking function with displacement $k$, every rearrangement of the parking function also has displacement $k$. For example, the parking functions $(1,1,2)$, $(1,2,1)$, and $(2,1,1)$ all have displacement~$2$. However, while the displacement is the same,
	these displacements arise in different ways: the displacement of $(1,1,2)$ arises from the sum $0+1+1=2$, while the displacements of $(1,2,1)$ and $(2,1,1)$ arise from the sum $0+0+2=2$.

	Classical parking functions with displacement 1 always arise from a single car parking one spot away from its preference. Aguillon, Alvarenga, Harris, Kotapati, Mart\'{i}nez Mori, Monroe, Saylor, Tieu, and Williams II established a bijection between classical parking functions with displacement $1$ and ideal states on the Tower of Hanoi game~\cite{Aguillon2022OnPF}. Unit-interval parking functions are classical parking functions in which all cars have displacement at most 1. Chaves Meyles, Harris, Jordaan, Rojas Kirby, Sehayek, and Spingarn established a bijection between unit-interval parking functions with $k$ cars having displacement $1$ and the set of $k$-dimensional facets of the permutohedron~\cite{unit_perm}. In this case, the displacement multiset is a multiset consisting of
	$k$ ones.

	The above-referenced work motivates the present study, where we introduce two new versions of parking functions inspired by the board games Monopoly and backgammon, in which cars finding their preferred spot occupied either return to the start of the street or bump the car occupying that spot to the start of the street, respectively. The car that returns to the start of the street then drives forward and parks in the first unoccupied spot it encounters. Next, we give the technical definitions of these parking schemes and present our major results related to the number of preference lists with specified parking outcomes, as well as our results related to displacement multisets.

	\subsection{Monopoly Parking Scheme}
	In the board game Monopoly~\cite{Monopoly}, there are Community Chest cards that instruct a player to ``Advance to Go'' which sends the player to the start of the board. Motivated by this feature of the game, we define the \textit{Monopoly parking scheme} as follows:
	\begin{itemize}
		\item Cars enter sequentially from car $1$ to $n$, in that order.

		\item Each car attempts to park in their preferred parking spot if available.

		\item If a car's preferred spot is occupied, then the car returns to the start of the street and parks in the first available spot it encounters.
	\end{itemize}
	If a preference list $\alpha=(a_{1},a_{2},\ldots,a_{n})\in[n]^{n}$ allows all cars to park under the Monopoly parking scheme, then $\alpha$ is a \textit{Monopoly parking function}. We denote the parking outcome of a preference list $\alpha$ under the Monopoly parking scheme by $\Out_{M}(\alpha)=\pi_{1}\pi_{2}\cdots\pi_{n}$, where $\pi_{j}=k$ denotes that car $k$ parks in the $j$-th spot. For succinctness, we refer to the outcome of a preference list $\alpha$ using the Monopoly parking scheme as the \textit{Monopoly outcome of~$\alpha$}. For example, using the preference list $\alpha=(1,4,4,3,2,2)$ of length $6$, we illustrate the Monopoly outcome of $\alpha$ in \Cref{fig:example of Monopoly pf}.

	\begin{figure}
		[h!]
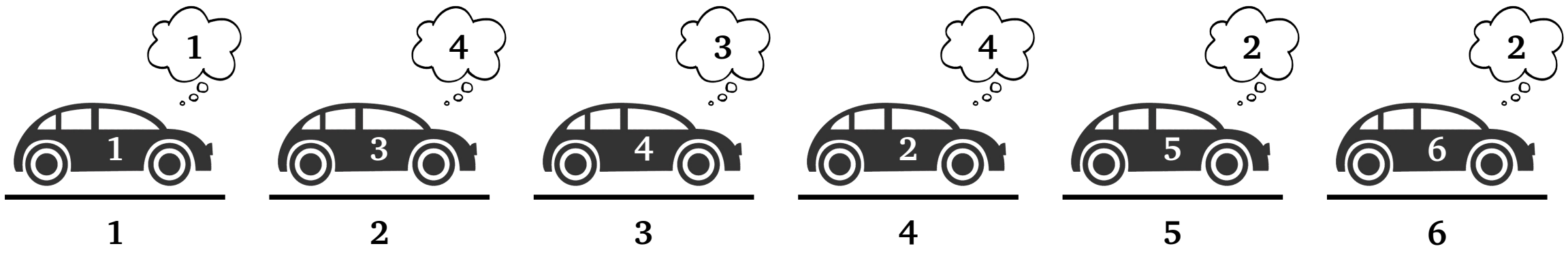
 \centering \resizebox{\textwidth}{!}{ \begin{tikzpicture}\node at(0,0){\includegraphics[width=1in]{car.png}}; \node at(3,0){\includegraphics[width=1in]{car.png}}; \node at(6,0){\includegraphics[width=1in]{car.png}}; \node at(9,0){\includegraphics[width=1in]{car.png}}; \node at(12,0){\includegraphics[width=1in]{car.png}}; \node at(15,0){\includegraphics[width=1in]{car.png}}; \node at(0,-.05){\textcolor{white}{\textbf{1}}}; \node at(3,-.05){\textcolor{white}{\textbf{3}}}; \node at(6,-.05){\textcolor{white}{\textbf{4}}}; \node at(9,-.05){\textcolor{white}{\textbf{2}}}; \node at(12,-.05){\textcolor{white}{\textbf{5}}}; \node at(15,-.05){\textcolor{white}{\textbf{6}}}; \node at(.9,1.){\includegraphics[width=.75in]{callout.png}}; \node at(3.9,1.){\includegraphics[width=.75in]{callout.png}}; \node at(6.9,1.){\includegraphics[width=.75in]{callout.png}}; \node at(9.9,1.){\includegraphics[width=.75in]{callout.png}}; \node at(12.9,1.){\includegraphics[width=.75in]{callout.png}}; \node at(15.9,1.){\includegraphics[width=.75in]{callout.png}}; \node at(.9,1.1){\textbf{1}}; \node at(3.9,1.1){\textbf{4}}; \node at(6.9,1.1){\textbf{3}}; \node at(9.9,1.1){\textbf{4}}; \node at(12.9,1.1){\textbf{2}}; \node at(15.9,1.1){\textbf{2}}; \draw[ultra thick](-1.25,-.6)--(1.25,-.6); \draw[ultra thick](1.75,-.6)--(4.25,-.6); \draw[ultra thick](4.75,-.6)--(7.25,-.6); \draw[ultra thick](7.75,-.6)--(10.25,-.6); \draw[ultra thick](10.75,-.6)--(13.25,-.6); \draw[ultra thick](13.75,-.6)--(16.25,-.6); \node at(0,-1){\textbf{1}}; \node at(3,-1){\textbf{2}}; \node at(6,-1){\textbf{3}}; \node at(9,-1){\textbf{4}}; \node at(12,-1){\textbf{5}}; \node at(15,-1){\textbf{6}};\end{tikzpicture} }
		\caption{The
		parking
		outcome
		for
		the
		preference
		list
		$(1,4,4,3,2,2)$
		under
		the
		Monopoly
		parking
		scheme.
		Hence,
		$\Out_{M}((1,4,4,3,2,2))=134256$.}
		\label{fig:example of Monopoly pf}
	\end{figure}

	We let $\MPF_{n}$ denote the set of all Monopoly parking functions. It is not difficult to see that \textit{every} preference list $\alpha\in[n]^{n}$ is a Monopoly parking function. We state this formally next, and prove it in \Cref{sec:Monopoly}.

	\begin{restatable}
		[]{proposition}{Monopoly} \label{prop: everything is Monopoly} If $n\in \NN$, then $\MPF_{n}=[n]^{n}$ and $|\MPF_{n}|=n^{n}$.
	\end{restatable}

	Given that every tuple in $[n]^{n}$ is a Monopoly parking function, we now ask: \textit{Given a permutation $\pi\in\Sym_{n}$, what tuples in $[n]^{n}$ have Monopoly outcome $\pi$?} We denote this set
	\[
		\Out_{M}^{-1}(\pi)\coloneqq \{\alpha\in[n]^{n}:\Out_{M}(\alpha)=\pi\}
	\]
	and as usual we let $|\Out_{M}^{-1}(\pi)|$ denote the cardinality of the set. To determine what is in the set $\Out_{M}^{-1}(\pi)$ and its cardinality, we introduce the following definition.
	\begin{definition}
		\label{def:goodindices} For an outcome $\pi=\pi_{1}\pi_{2}\cdots\pi_{n}\in\Sym_{n}$, the \textit{set of good indices of $\pi$} is defined as the set $I_{\pi}\subseteq[n]$ satisfying the following condition: $i\in I_{\pi}$ if and only if $\pi_{j}<\pi_{i}$ for all $j<i$.
	\end{definition}

	In other words, the set of good indices of an outcome consists of all parking spots $i$ such that the car in spot $i$ is numbered higher than every car parked to its left. Note that $1$ is vacuously in the set of good indices for every permutation $\pi$, as there are no indices $j<1$. As an example, note that the set of good indices for $\pi=532618794$ is $I_{\pi}=\{1,4,6,8\}$. Then, we have the following result for Monopoly parking functions.

	\begin{restatable}
		[]{theorem}{Monopolyoutcomes} \label{thm: everything is Monopoly} Fix an outcome $\pi=\pi_{1}\pi_{2}\cdots \pi_{n}$ of length $n$ and let $I_{\pi}$ be the set of good indices of $\pi$. Then the number of preference lists $\alpha=(a_{1},a_{2},\dots, a_{n})\in[n]^{n}$ whose Monopoly outcome is $\pi$ is given by
		\[
			|\Out_{M}^{-1}(\pi)|=\displaystyle\prod_{i\in I_\pi}\pi_{i}.
		\]
	\end{restatable}
	By \Cref{thm: everything is Monopoly}, given that the set of good indices of $\pi=532618794$ is $I_{\pi}=\{1,4,6,8\}$, we have that the number of preference lists in $[n]^{n}$ whose Monopoly outcome is $\pi$ is given by $\pi_{1}\cdot\pi_{4}\cdot\pi_{6}\cdot\pi_{8}={ 5\cdot 6\cdot 8\cdot 9}=2160$.

	The key idea in the proof of \Cref{thm: everything is Monopoly} (see Section \ref{sec:Monopoly}) is that each car $k$ contributes either $1$ or $k$ to the product. Consequently, the maximum possible product is $n!$, achieved by the outcome $12\cdots(n-1)n$, while the minimum possible product is $n$, achieved by any outcome $\pi=\pi_{1}\pi_{2}\cdots\pi_{n}$ with $\pi_{1}=n$.

	Thus, we readily arrive at the following enumerative results:

	\begin{corollary}
		\label{lem:increase and decrease outcomes in Monopoly} Let $n\geq 1$.
		\begin{itemize}
			\item If $\pi=123\cdots (n-1)n$ is the identity permutation, then $|\Out_{M}^{-1}(123\cdots (n-1)n)|=n!$.

			\item If $\pi=\pi_{1}\pi_{2}\dots\pi_{n}$ is any outcome such that $\pi_{1}=n,$ then $|\Out_{M}^{-1}\left(\pi\right)|=n.$ In particular, $|\Out_{M}^{-1}\left(n (n-1)\cdots 321\right)|=n$.
		\end{itemize}
	\end{corollary}
    
	\subsection{Backgammon Parking Scheme}
	In the game of backgammon~\cite{Backgammon}, checkers are placed on a board. When a checker lands on a space occupied by exactly one of the opponent’s checkers, it sends the opponent's checker off the board. On the opponent's next turn, they get a chance to re-enter their checker onto the board.

	Inspired by this feature of backgammon, we introduce the \textit{backgammon parking scheme}:
	\begin{itemize}
		\item Cars enter sequentially from car $1$ to $n$, in that order.

		\item Cars attempt to park in their preferred parking spot if available.

		\item If a car's preferred spot is occupied, then the car bumps the parked car out of its preferred parking spot. The bumped car returns to the start of the street and parks in the first available spot it encounters.
	\end{itemize}
	If the preference list $\alpha=(a_{1},a_{2},\ldots,a_{n})\in[n]^{n}$ allows all cars to park under the backgammon parking scheme, then $\alpha$ is a \textit{backgammon parking function of length $n$}. We denote the parking outcome of a preference list $\alpha$ under the backgammon parking scheme by $\Out_{B}(\alpha)=\pi_{1}\pi_{2}\cdots\pi_{n}$, where $\pi_{j}=k$ denotes that car $k$ parks in the $j$-th spot. For succinctness, we refer to the outcome of a preference list $\alpha$ using the backgammon parking scheme as the \textit{backgammon outcome of $\alpha$}. For example, using the preference list $\alpha=(1,4,4,3,2,2)$ of length 6, we illustrate the backgammon outcome of $\alpha$ in \Cref{fig:backgammon example}. We let $\BPF_{n}$ denote the set of backgammon parking functions of length
	$n$.
	\begin{figure}
		[h!]
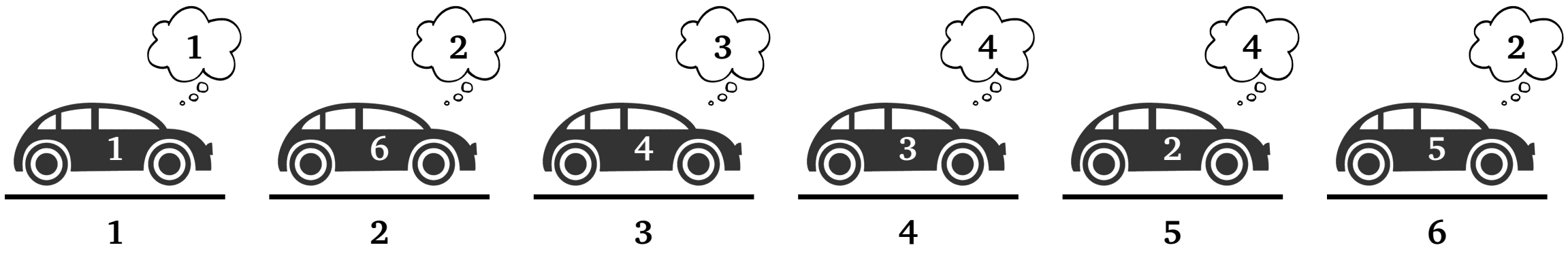
 \centering

		\resizebox{\textwidth}{!}{ \begin{tikzpicture}\node at(0,0){\includegraphics[width=1in]{car.png}}; \node at(3,0){\includegraphics[width=1in]{car.png}}; \node at(6,0){\includegraphics[width=1in]{car.png}}; \node at(9,0){\includegraphics[width=1in]{car.png}}; \node at(12,0){\includegraphics[width=1in]{car.png}}; \node at(15,0){\includegraphics[width=1in]{car.png}}; \node at(0,-.05){\textcolor{white}{\textbf{1}}}; \node at(3,-.05){\textcolor{white}{\textbf{6}}}; \node at(6,-.05){\textcolor{white}{\textbf{4}}}; \node at(9,-.05){\textcolor{white}{\textbf{3}}}; \node at(12,-.05){\textcolor{white}{\textbf{2}}}; \node at(15,-.05){\textcolor{white}{\textbf{5}}}; \node at(.9,1.){\includegraphics[width=.75in]{callout.png}}; \node at(3.9,1.){\includegraphics[width=.75in]{callout.png}}; \node at(6.9,1.){\includegraphics[width=.75in]{callout.png}}; \node at(9.9,1.){\includegraphics[width=.75in]{callout.png}}; \node at(12.9,1.){\includegraphics[width=.75in]{callout.png}}; \node at(15.9,1.){\includegraphics[width=.75in]{callout.png}}; \node at(.9,1.1){\textbf{1}}; \node at(3.9,1.1){\textbf{2}}; \node at(6.9,1.1){\textbf{3}}; \node at(9.9,1.1){\textbf{4}}; \node at(12.9,1.1){\textbf{4}}; \node at(15.9,1.1){\textbf{2}}; \draw[ultra thick](-1.25,-.6)--(1.25,-.6); \draw[ultra thick](1.75,-.6)--(4.25,-.6); \draw[ultra thick](4.75,-.6)--(7.25,-.6); \draw[ultra thick](7.75,-.6)--(10.25,-.6); \draw[ultra thick](10.75,-.6)--(13.25,-.6); \draw[ultra thick](13.75,-.6)--(16.25,-.6); \node at(0,-1){\textbf{1}}; \node at(3,-1){\textbf{2}}; \node at(6,-1){\textbf{3}}; \node at(9,-1){\textbf{4}}; \node at(12,-1){\textbf{5}}; \node at(15,-1){\textbf{6}};\end{tikzpicture} }
		\caption{ The parking outcome for the preference list $(1,4,4,3,2,2)$ under the backgammon parking scheme. Hence, $\Out_{B}((1,4,4,3,2,2))=164325$. }
		\label{fig:backgammon example}
	\end{figure}

	Although similar, the definition of backgammon parking functions is different from that of Most Valuable Player (MVP) parking functions introduced by Harris, Kamau, Mart\'inez Mori, and Tian~\cite{MVP}. Under the MVP parking scheme, when a car finds its preferred spot occupied, it bumps the current occupant out of that spot. The bumped car then seeks forward on the street for a new parking spot, and parking may fail if no spots are available past its original preference. If a preference list parks all of the cars under the MVP parking scheme, we say that it is an \textit{MVP parking function}. In contrast, under the backgammon parking scheme, a bumped car returns to the start of the street and parks in the first available spot it encounters, which could be earlier than its preference.

	The number of MVP parking functions of length $n$ is $(n+1)^{n-1}$ (see~\cite[Corollary 2.1]{MVP}), which is equal to the number of classical parking functions. Our first result for backgammon parking functions establishes that every preference list in $[n]^{n}$ is a backgammon parking function of length $n$; we prove this in \Cref{sec:Backgammon}.

	\begin{restatable}
		[]{proposition}{backgammon} \label{prop: everything is backgammon} If $n\in\NN$, then $\BPF_{n}=[n]^{n}$ and $|\BPF_{n}|=n^{n}$.
	\end{restatable}

	As in the Monopoly case, we are interested in determining the preference lists with a fixed backgammon outcome. For $\pi\in\Sym_{n}$, we denote this set by
	\[
		\Out_{B}^{-1}(\pi)\coloneqq \{\alpha\in[n]^{n}:\Out_{B}(\alpha)=\pi\}
	\]
	and as usual we let $|\Out_{B}^{-1}(\pi)|$ denote the cardinality of the set.

	For Monopoly parking functions, we were able to fully enumerate the preference lists resulting in a specified parking outcome. However, the backgammon parking scheme returning a previously parked car to the start of the street creates some challenges in addressing the analogous question. The complication arising from cars being bumped from their spot was also a major challenge in the MVP parking function problem, where the authors establish an upper bound for the number of preference lists that park the cars in a given order~\cite[Theorem 3.1]{MVP}.

	We prove that for any permutation $\pi$, the number of preference lists with backgammon outcome equal to $\pi$ is at most $C_{n}=\frac{1}{n+1}\binom{2n}{n}$, the $n$-th Catalan number, with corresponding OEIS sequence \seqnum{A000108}. Furthermore, we determine the outcomes for which this upper bound is an equality and for which it is a strict inequality.

	\begin{restatable}
		[]{theorem}{AtMostCatalan} \label{thm:at most Catalan} Fix $n\in \NN$ and let $C_{n}$ denote the $n$-th Catalan number. For any $\pi\in\Sym_{n}$, we have
		\[
			|\Out_{B}^{-1}(\pi)|\leq C_{n}.
		\]
		Furthermore, $\left|\Out_{B}^{-1}\left(\pi\right)\right|=C_{n}$ whenever $\pi=12\cdots (n), 12\cdots(n-2)(n)(n-1)$, or $12\cdots(n-3)(n)(n-2)(n-1)$. Otherwise, the inequality is strict.
	\end{restatable}

	We also establish that for any permutation $\pi$, there are at least $2^{n-1}$ preference lists in $[n]^{n}$ whose backgammon outcome is equal to $\pi$, thereby providing a lower bound for this count.

	\begin{restatable}
		{theorem}{BackgammonLowerBound} \label{thm:backgammon lower bound} For any $\pi\in\Sym_{n}$, we have $|\Out_{B}^{-1}(\pi)|\geq 2^{n-1}$. For $n\geq 5$, this inequality is strict for all $\pi\in \Sym_{n}$.
	\end{restatable}
	We give proofs of \Cref{thm:at most Catalan} and \Cref{thm:backgammon lower bound} in \Cref{subsec:bounds for backgammon}.

	As noted above, it is difficult to determine the number of backgammon parking functions with a fixed outcome in general, but for some special cases we establish connections to the Motzkin numbers, which enumerate Motzkin lattice paths (\Cref{def:Motzkin paths}). Here $M_{n}$ denotes the $n$-th Motzkin number (see \Cref{eq:motzkin}) and $\{M_{n}\}_{n\geq 1}$ corresponds to OEIS sequence \seqnum{A001006}. In \Cref{sec:dec perms in backgammon}, we establish the following results.

	\begin{itemize}[leftmargin=0in]
		\item[]
			\begin{restatable}
				{theorem}{thmMotzkin} \label{thm:motzkin} For $n\ge 1$, $\left|\Out_{B}^{-1}\left(n(n-1)\cdots21\right)\right|=M_{n}$.
			\end{restatable}

		\item[]
			\begin{restatable}
				[]{corollary}{differenceofadjacentmotzkin} \label{cor:difference of adjacent motzkin} For $n\ge 1$, $\left|\Out_{B}^{-1}\left(1(n)(n-1)\cdots32\right)\right|=M_{n+1}-M_{n}$.
			\end{restatable}

		\item[]
			\begin{restatable}
				[]{corollary}{moremotzkin} \label{cor:more motzkin} For $n\ge 2$, $\left|\Out_{B}^{-1}\left(12(n)(n-1)\cdots 3\right)\right|=M_{n+2}-2M_{n+1}+M_{n}-M_{n-1}.$
			\end{restatable}
	\end{itemize}

	\subsection{Displacement}

	Recall that the displacement $d_{i}$ of car $i$ is the \textit{absolute value} of the difference between the spot car $i$ prefers and where it ultimately parks. For any parking scheme, we say a car is ``lucky'' if at the conclusion of parking, it ends up in the spot that it preferred. Otherwise, we say that the car is ``unlucky.'' Thus, the unlucky cars are precisely those with positive displacement.

	Let $X\in\{C,M,B\}$ be a parking scheme, where $X=C$ denotes the classical parking scheme, $X=M$ denotes the Monopoly parking scheme, and $X=B$ denotes the backgammon parking scheme. Then we define:
	\begin{definition}
		\label{def:displacementvector} For a preference list $\alpha=(a_{1},a_{2},\dots, a_{n})\in [n]^{n}$ with outcome $\Out_{X}(\alpha)=\pi_{1}\pi_{2}\cdots\pi_{n}$ and $k\geq 0$ unlucky cars, the \textit{displacement multiset} $\calD_{X}(\alpha)$ is the multiset of the displacements $d_{i}$ of each unlucky car $i$.
	\end{definition}

	To highlight how \Cref{def:displacementvector} depends on the parking scheme used, we present some examples. Given an outcome $\pi=\pi_{1}\pi_{2}\cdots \pi_{n}\in\Sym_{n}$, we define $\pi^{-1}=p_{1}p_{2}\cdots p_{n}$ to be the permutation of $[n]$ such that $\pi_{p_j}=j$ for all $j\in [n]$; that is, $p_{j}=k$ denotes that car $j$ parks in spot $k$.
	\begin{example}
		\label{ex:displacement}As illustrated in \Cref{fig:example of pf,fig:example of Monopoly pf,fig:backgammon example}, the preference list $\alpha=(1,4,4,3,2,2)$ satisfies the following:
		\begin{itemize}[leftmargin=.25in]
			\item The (classical) outcome of $\alpha$ is $\Out_{C}(\alpha)=\pi=154236$, so $\pi^{-1}=145326$ and the (classical) displacement multiset is $\calD_{C}(\alpha)=\{1,4\}$.

			\item The Monopoly outcome of $\alpha$ is $\Out_{M}(\alpha)=\pi=134256$, so $\pi^{-1}=142356$ and the Monopoly displacement multiset is $\calD_{M}(\alpha)=\{2,3,4\}$.

			\item The backgammon outcome of $\alpha$ is $\Out_{B}(\alpha)=\pi=164325$, so $\pi^{-1}=154362$ and the backgammon displacement multiset is $\calD_{B}(\alpha)=\{1,4\}$.
		\end{itemize}
	\end{example}

	\begin{example}
		\label{ex:displacement2} Similarly, for the preference list $\alpha=(1,3,3,1),$ we have:
		\begin{itemize}[leftmargin=.25in]
			\item The (classical) outcome of $\alpha$ is $\Out_{C}(\alpha)=\pi=1423$, so $\pi^{-1}=1342$ and the (classical) displacement multiset is $\calD_{C}(\alpha)=\{1,1\}$.

			\item The Monopoly outcome of $\alpha$ is $\Out_{M}(\alpha)=\pi=1324$, so $\pi^{-1}=\pi$ and the Monopoly displacement multiset is $\calD_{M}(\alpha)=\{1,3\}$.

			\item The backgammon outcome of $\alpha$ is $\Out_{B}(\alpha)=\pi=4231$, so $\pi^{-1}=\pi$ and the backgammon displacement multiset is $\calD_{B}(\alpha)=\{1,3\}$.
		\end{itemize}
	\end{example}

	We define $N_{X}(\calD_{X})$ to be the number of parking functions with displacement multiset $\calD_{X}$ using the parking scheme $X\in\{C,M,B\}$. For the \textbf{Monopoly parking scheme}, we prove the following result. We use $H_k$ to denote the $k$-{th} Harmonic number $\sum_{j=1}^k\frac{1}{j}$.
	\begin{restatable}
		[]{theorem}{Monopolyreluctant}\label{thm:Monopoly1reluctant} For $x\in [n-1]$, the number of Monopoly parking functions $\alpha\in\MPF_{n}$ with displacement multiset $\calD_{M}(\alpha)=\{x\}$ is \[N_M(\{x\})=H_{n}-H_{x}+H_{n-x+1}-1.\]
	\end{restatable}
	
    For the \textbf{backgammon parking scheme}, we prove the following results.

	\begin{restatable}
		[]{theorem}{backgammondisplacement} \label{thm:backgammon displacement 1} For $x\in [n-1]$, the number of backgammon parking functions $\alpha\in\BPF_{n}$ with displacement multiset $\calD_{B}(\alpha)=\{x\}$ is
		\[
			N_{B}(\{x\})=(n-x)n!.
		\]
	\end{restatable}

	\begin{restatable}
		[]{theorem}{backgammondisplacementtwo} For $x_{1}, x_{2}\in[n-1]$ with $x_{1}> x_{2}$, the number of backgammon parking functions $\alpha\in\BPF_{n}$ with displacement multiset $\calD_{B}(\alpha)=\{x_{1}, x_{2}\}$ is
		\begin{align}
			\label{eq:x1 neq x2}\raisetag{5ex}\displaystyle N_{B}{(\{x_1,x_2\})} & = \begin{cases}\dfrac{n!}{6}\left(3n^2-(3x_1+3x_2+8)n+(3x_1x_2+8x_1+4x_2)\right) & \text{ when } n\ge x_1+x_2,\\[1.2em] \dfrac{n!}{6}\left(3n^2-(3x_1+3x_2+4)n+(3x_1x_2+4x_1)\right) & \text{ otherwise. }\end{cases}
		\end{align}

		For $x_{1}\in [n-1]$, the number of backgammon parking functions $\alpha\in\BPF_{n}$ with displacement multiset $\calD_{B}(\alpha)=\{x_{1},x_{1}\}$ is
		\begin{align}
			\label{eq:x1 eq x2}N_{B}{(\{x_1,x_1\})}= \begin{cases}\dfrac{n!}{12}\left(3n^2-(6x_1+7)n+(3x_1^2+11x_1)\right)& \text{ when } n\ge 2x_1,\\ &\\ \dfrac{n!}{4}\left(n^2-(2x_1+1)n+(x_1^2+x_1)\right) & \text{ otherwise. }\end{cases}
		\end{align}
	\end{restatable}

	We also prove the following result comparing the number of backgammon and Monopoly parking functions for the same displacement multiset of size 1.
	\begin{restatable}
		[]{theorem}{displacementcomparison}\label{thm:displacementcomparison} For $x\in [n-1],$ the number $N_{M}(\{x\})$ of Monopoly parking functions in $\MPF_{n}$ with displacement multiset $\{x\}$ is at most the number $N_{B}(\{x\})$ of backgammon parking functions in $\BPF_{n}$ with the displacement multiset $\{x\}$. That is, we have,
		\[
			N_{M}(\{x\})\leq N_{B}(\{x\}).
		\]
		Moreover, equality holds if and only if $n=2$ and $x=1$.
	\end{restatable}
	\subsection{Organization}
	In the sections that follow, we consider each of the newly introduced parking schemes and prove our major results. \Cref{sec:Monopoly} contains our results for counting the number of Monopoly parking functions with a specific outcome, and \Cref{sec:Backgammon} provides analogous results for backgammon parking functions. In \Cref{sec:displacement}, we give results on displacement multisets for both Monopoly and backgammon parking functions. We conclude with \Cref{sec:future}, where we detail some directions for future work.

	\begin{remark}
		Readers interested primarily in Monopoly parking functions may consult Sections \ref{sec:Monopoly} and \ref{sub:dispMon}, together with the comparison in Section \ref{sub:dispComp}. Readers interested primarily in backgammon parking functions may consult Sections \ref{sec:Backgammon}, \ref{sub:dispBack}, and \ref{sub:dispComp}. Section \ref{sec:future} discusses open problems and connections with MVP parking functions.
	\end{remark}

	\section{Monopoly Parking Functions}
	\label{sec:Monopoly} We begin by giving a proof that every tuple in $[n]^{n}$ is a Monopoly parking function.

	\Monopoly*
	\begin{proof}
		The set containment $\MPF_{n}\subseteq[n]^{n}$ holds by definition. Now consider an arbitrary tuple $(a_{1},a_{2},\ldots,a_{n})\in[n]^{n}$. We claim that the tuple allows all cars to park using the Monopoly parking scheme. The first car always parks in its preferred spot. Suppose that at step $i\in [n],$ there are $i-1$ cars parked in $i-1$ spots, and there are $n-i+1$ open spots. Once car $i$ enters, it either parks in its preferred spot or is forced to return to the start of the street and parks in the first available spot it encounters. Since $n-i+1>0$, such a spot always exists. Thus, after $n$ steps, all cars end up parked in the $n$ spots.
	\end{proof}

	Recall that for any outcome $\pi$, the set of good indices $I_{\pi}$ is given by $\{\,i\in[n]\colon \pi_{j}< \pi_{i}\text{ for all }j<i\,\}$ (see \Cref{def:goodindices}). Thus, we prove \Cref{thm: everything is Monopoly}. \Monopolyoutcomes*

	\begin{proof}
		We show that for each good index $i\in I_{\pi}$, car $\pi_{i}$ contributes a factor of $\pi_{i}$ to the count, while if $i$ is not a good index, then car $\pi_{i}$ contributes a factor of $1$ to the count.

		Suppose $i$ is not a good index. Then there exists some $j<i$ such that $\pi_{j}>\pi_{i}$. Since car $\pi_{j}$ arrives after car $\pi_{i}$, when car $\pi_{i}$ arrives, it finds at least one spot earlier than spot $i$ that is vacant, namely spot $j<i$. Therefore, in order for car $\pi_{i}$ to end up in spot $i$, it must have preferred spot $i$ to begin with, thus contributing a factor of $1$ to the count. On the other hand, if $i$ is a good index, then $\pi_{j}<\pi_{i}$ for every $j<i$, which means that all cars that eventually occupy spots $1,\dots,i-1$ arrive before car $\pi_{i}$. Hence, when car $\pi_{i}$ arrives, every spot before $i$ is already occupied and spot $i$ is available. Therefore, car $\pi_{i}$ will park in spot $i$ whenever its preferred spot is either $i$ itself or
		any of the $\pi_{i}-1$ previously occupied spots. Thus, car $\pi_{i}$ contributes a factor of $\pi_{i}$, and the total count is given by the product $\displaystyle\prod_{i\in I_\pi}\pi_{i}$ where $I_{\pi}=\{\,i\in\ [n] : \pi_{j}< \pi_{i}\text{ for all }j<i\,\}$ is the set of good indices of $\pi$ (Definition \ref{def:goodindices}). This proves the claim. \qedhere
	\end{proof}

	\section{Backgammon Parking Functions}
	\label{sec:Backgammon} We begin by giving a proof of \Cref{prop: everything is backgammon}.

	\backgammon*

	\begin{proof}
		The first car always parks in its preferred spot. Suppose that at step $i\in [n],$ there are $i-1$ cars parked in $i-1$ spots, and there are $n-i+1$ open spots. Once car $i$ enters, it either finds its preferred spot unoccupied (and parks) or bumps a previously parked car $j<i$ out of car~$j$'s current spot. In the second case, since car $i$ has taken up an already occupied spot, there are still $n-i+1>0$ open spots, all of which are available for car~$j$ to take. As car $j$ re-enters the street at the beginning, it is guaranteed to find an empty spot in which to park. Thus, after $n$ steps, all cars will be parked in the $n$ spots.
	\end{proof}

	\subsection{Bounds for the Number of Backgammon Parking Functions with a Fixed Outcome}
	\label{subsec:bounds for backgammon} Similar to Theorem \ref{thm: everything is Monopoly} for Monopoly parking functions, we wish to answer the following question: given an outcome $\pi=\pi_{1}\pi_{2}\cdots \pi_{n}\in\Sym_{n}$, how many preference lists $\alpha=(a_{1},a_{2},\dots, a_{n})\in[n]^{n}$ allow cars to park in this configuration using the {{backgammon parking scheme}}? Theorem \ref{thm:at most Catalan} gives an upper bound for this count. Before restating this result, we recall some background on Dyck paths.
	\begin{definition}
		A \textit{Dyck path} of semilength $n$ (or length $2n$) is a lattice path from $(0,0)$ to $(n,n)$ that lies on or above the main diagonal $y=x$ and consists of East steps (denoted by $\E$) and North steps (denoted by $\N$). We let $D_{n}$ denote the set of Dyck paths of semilength $n$.
	\end{definition}
	In \Cref{fig:Dyck paths n=3}, we illustrate the five Dyck paths of semilength $3$. We recall that Dyck paths of semilength $n$ are Catalan objects and, as such, $|D_{n}|=C_{n}=\frac{1}{n+1}\binom{2n}{n}$, where $C_{n}$ denotes the $n$-th Catalan number. For more on Catalan objects, we recommend Stanley's book on Catalan numbers~\cite{StanleyCatalan}. Moreover, we recall that there is a bijection between parking functions and Dyck paths with north steps labeled by the set $[n]$ in such a way that consecutive north steps are labeled in increasing order. This bijection is described in~\cite[Theorem 6]{bib:DescentsInPF}.
	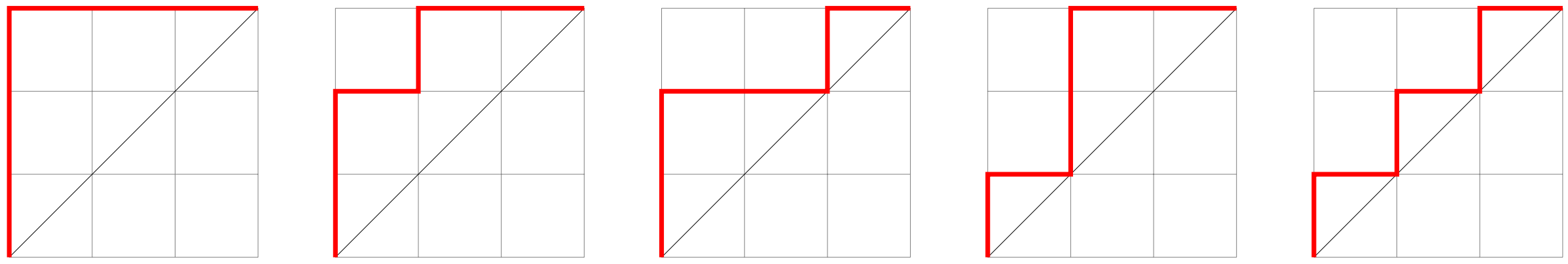
\begin{figure}
		\centering \resizebox{\textwidth}{!}{ \begin{tikzcd}\draw[step=1cm,gray,very thin] (0,0) grid (3,3); \draw[black, very thin] (0,0)--(3,3); \draw[ultra thick, red] (0,0)--(0,3)--(3,3); \node at(1.5,-.5){\large{\N\N\N\E\E\E}};\end{tikzcd} \qquad \begin{tikzcd}\draw[step=1cm,gray,very thin] (0,0) grid (3,3); \draw[black, very thin] (0,0)--(3,3); \draw[ultra thick, red] (0,0)--(0,2)--(1,2)--(1,3)--(3,3); \node at(1.5,-.5){\large{\N\N\E\N\E\E}};\end{tikzcd} \qquad \begin{tikzcd}\draw[step=1cm,gray,very thin] (0,0) grid (3,3); \draw[black, very thin] (0,0)--(3,3); \draw[ultra thick, red] (0,0)--(0,2)--(2,2)--(2,3)--(3,3); \node at(1.5,-.5){\large{\N\N\E\E\N\E}};\end{tikzcd} \qquad \begin{tikzcd}\draw[step=1cm,gray,very thin] (0,0) grid (3,3); \draw[black, very thin] (0,0)--(3,3); \draw[ultra thick, red] (0,0)--(0,1)--(1,1)--(1,3)--(3,3); \node at(1.5,-.5){\large{\N\E\N\N\E\E}};\end{tikzcd} \qquad \begin{tikzcd}\draw[step=1cm,gray,very thin] (0,0) grid (3,3); \draw[black, very thin] (0,0)--(3,3); \draw[ultra thick, red] (0,0)--(0,1)--(1,1)--(1,2)--(2,2)--(2,3)--(3,3); \node at(1.5,-.5){\large{\N\E\N\E\N\E}};\end{tikzcd} }
		\caption{Illustration of the five distinct Dyck paths of semilength $3$.}
		\label{fig:Dyck paths n=3}
	\end{figure}

	Another classical set of combinatorial objects enumerated by the Catalan numbers is the possible arrangements of $n$ pairs of parentheses that are valid. A bijection between valid placements of $n$ pairs of parentheses and a Dyck path of semilength $n$ is given by replacing every $\N$ with ``$($'' and $\E$ with ``$)$''. For example, the Dyck path of semilength $3$ given by $\N\E\N\N\E\E$ corresponds to the valid parentheses ()(()). Henceforth, we use the term ``Dyck path'' to refer to both the lattice path and its corresponding valid parenthesization.

	\begin{definition}
		For each preference list $\alpha$ for $n$ cars, we associate to it a set partition of $[n]$ known as the \textit{content partition} of $\alpha$ and denoted by $\calS(\alpha)$, where the blocks of $\calS(\alpha)$ correspond to maximal subsets of cars that share the same preference.
	\end{definition}
	\begin{example}
		If $\alpha=(1,2,2,4,4,3)$, then $\calS(\alpha)=\{\{1\},\{2,3\},\{4,5\},\{6\}\}$.
	\end{example}
	\begin{lemma}
		\label{lem:recoveralpha} Given a set partition $S$ of $[n]$ and a permutation $\pi$ of $[n]$, there is at most one $\alpha\in[n]^{n}$ such that $S=\calS(\alpha),$ and $\pi=\Out_{B}(\alpha)$.
	\end{lemma}
	\begin{proof}
		If no such preference list exists, we are done, so assume at least one exists and call it $\alpha$. Recall that $\Out_{B}(\alpha)=\pi_{1}\pi_{2}\cdots\pi_{n}$, where $\pi_{x}=y$ means that car $y$ parked in spot $x$. Let $J$ be a block of $\calS(\alpha)$, and let $k=\max(J)$ be the highest index in $J$. Then, car $k$ parks in its preferred position, and every car indexed by a number in $J$ shares its preference. Since car $k$ parks in spot $\pi^{-1}(k)$, this common preference is $\pi^{-1}(k)$. As the block $J$ was arbitrary, we have that for each block $J'\in \calS(\alpha)$, all cars indexed by elements in block $J'$ share the preference $\pi^{-1}(\max(J'))$. Thus, $\alpha$ is completely determined by $\calS(\alpha)=S$ and $\Out_{B}(\alpha)=\pi$, which finishes the proof.
	\end{proof}

	\begin{example}
		Consider $S=\{\{1\},\{2,3\},\{4,5\},\{6\}\}$ and $\pi=654213$. Using \Cref{lem:recoveralpha}, we can uniquely find $\alpha$ such that $\mathcal{S}(\alpha)=S$ and $\Out_{B}(\alpha)=\pi$, as follows: Since car $1$ parked in spot $5$ and it is the largest-numbered car in its own block, car $1$ has to have preferred spot $5$. Since car $3$ parks in spot $6$, and it is the largest-numbered car in its block, all cars in the same block (i.e., cars 2 and 3) have to prefer spot $6$. Since car $5$ parks in spot $2$, and it is the largest-numbered car in its block, all cars in the same block (i.e., cars 4 and 5) have to prefer spot $2$. Since car $6$ parks in spot $1$, and it is the largest-numbered car in its block, it has to prefer spot $1$. This then establishes that the only possible $\alpha$
		is $\alpha=(5,6,6,2,2,1)$. Indeed, we can check that this $\alpha$ parks the cars in the order $\pi=654213$.
	\end{example}

	\begin{example}
		Consider $S=\{\{1\},\{2,3\},\{4,5\},\{6\}\}$ and $\pi=652413$. The cars corresponding to the highest index in each block, namely cars $1,3,5,6$, must prefer the spots they park in. Car $2$ has the same preference as car $3$ and car $4$ has the same preference as car $5$, so the only potential preference list is again $\alpha=(5,6,6,2,2,1)$. However, this yields the outcome $654213$, which is not $\pi=652413$. Consequently, there is no preference list $\alpha$ such that $\mathcal{S}(\alpha)=S$ and $\Out_{B}(\alpha)=\pi$.
	\end{example}

	\AtMostCatalan*

	For clarity, we separate the proof of \Cref{thm:at most Catalan} into a few technical results: \Cref{prop:justupper}, \Cref{thm:strictupper}, and \Cref{lem:combined}.

	\begin{proposition}
		\label{prop:justupper} For any $\pi\in\Sym_{n}$, $\left|\Out_{B}^{-1}\left(\pi\right)\right|\le C_{n}$.
	\end{proposition}

	Before presenting the proof of this result, we introduce the notion of a partial parking outcome.
	\begin{definition}
		\label{def:partial parking outcome} Let $\alpha$ be a backgammon parking function with outcome $\Out_{B}(\alpha)=\pi$. Fix $t\in[0,n]$. After $t$ cars have parked on the street with preference list $\alpha$, we denote by $\Out_{B,t}(\alpha)=\tau_{1}\tau_{2}\cdots\tau_{n}$ the \textit{partial backgammon outcome after step $t$}, where
		\[
			\tau_{i}=
			\begin{cases}
				0&\mbox{if there is no car parked in spot $i$}\\ j&\mbox{if car $j$ is parked in spot $i$}.
			\end{cases}
		\]
	\end{definition}
	Note that for every $t\in[0,n]$, $\Out_{B,t}(\alpha)$ is a permutation of the multiset $\{\underbrace{0,\ldots,0}_{n-t},1,2,\ldots,t\}$. In particular, $\Out_{B,n}(\alpha)=\Out_{B}(\alpha)=\pi$ and $\Out_{B,0}(\alpha)=\underbrace{0\cdots0}_{n}.$ For the example in \Cref{fig:backgammon example}, the partial backgammon outcome of $\alpha=(1,4,4,3,2,2)$ after step $4$ is $124300$.

	\begin{remark}
		\label{rem:same def for partial parking outcome under Monopoly} Note that we can analogously define a partial parking outcome $\Out_{M,t}(\alpha)$ for the Monopoly parking scheme, though we will refer to this definition only in \Cref{sec:future}.
	\end{remark}

	Inspired by the discussion in~\cite{ideaforproof}, we introduce the notion of a Dyck path $P_{S}$ associated to a set partition $S$ of $[n]$, which we call the \textit{pathification} of $S.$

	\begin{definition}
		\label{def:dyckpath} Any set partition $S$ of $[n]$ can be mapped to a Dyck path $P_{S}$ of semilength $n$ in the following manner:
		\begin{itemize}
			\item If $\{i\}$ is a block of $S$, then the $(2i-1)$-th and $(2i)$-th steps of the Dyck path are $()$,

			\item If $i$ is in a block of $S$ of size greater than one and
				\begin{itemize}
					\item if $i$ is the smallest element in its block, then the $(2i-1)$-th and $(2i)$-th steps of the Dyck path are $(($,

					\item if $i$ is the largest element in its block, then the $(2i-1)$-th and $(2i)$-th steps of the Dyck path are $))$,

					\item if $i$ is neither the smallest nor largest element in its block, then the $(2i-1)$-th and $(2i)$-th steps of the Dyck path are $)($.
				\end{itemize}
		\end{itemize}
	\end{definition}
	Note that each block contributes an equal number of $($ and $)$. It is easy to verify that the output of this map is always a Dyck path. However, this map is not a bijection in general.

	For examples of the map, consider $S=\{\{1,2\},\{3\}\}$ and $S'=\{\{1,3,6\},\{2,5\},\{4\}\}$. Then $P_{S}=(())()$ and $P_{S'}=(((()(()))))$.

	We now state the following lemma.
	\begin{lemma}
		\label{lem:onePrefListperP} For every pair of a permutation $\pi\in\Sym_{n}$ and a Dyck path $P$ of semilength $n$ (i.e., length $2n$), there is at most one preference list $\alpha$ with outcome $\pi$ whose content partition $\mathcal{S(\alpha)}$ maps to $P$ under the pathification map (see \Cref{def:dyckpath}).
	\end{lemma}
	\begin{proof}
		By Lemma~\ref{lem:recoveralpha}, at most one preference list $\alpha$ with a given content partition can result in a particular outcome $\pi$. Therefore, it suffices to show that for every pair $(\pi, P)$, with $\pi$ a permutation and $P$ a Dyck path, there is at most one set partition $S$ mapping to $P$, such that $S=\calS(\alpha)$ for some $\alpha$ with $\Out_{B}(\alpha)=\pi$.

		To this end, fix a permutation $\pi$ and fix a Dyck path $P$. If $i$ corresponds to $))$, that is, the $(2i-1)$-th and $(2i)$-th steps of the Dyck path $P$ are $))$, then $i$ must be the last element of a block of $S$, while if the $(2i-1)$-th and $(2i)$-th steps of the Dyck path are $()$, then $\{i\}$ must be a block of size $1$. At this point, each block has had a representative chosen, so no new blocks will be added. We now proceed by repeatedly working backward from $n$ until every element has been assigned to a block. At each step, we consider the highest element, $j$, of $[n]$ which is not yet assigned to a block, so in particular, proceeding analogously to the proof of \Cref{lem:recoveralpha}, we uniquely know the preferences of the cars numbered higher than $j$. Next, we show
		that there is always at most one choice for which block to assign $j$ to, regardless of whether the $(2j-1)$-th and $(2j)$-th steps of $P$ are $(($ or $)($, i.e., whether $j$ is the minimal element of the block it is contained in or whether $j$ is some element that is neither the minimum nor maximum of the block it is contained in.

		Suppose, for the sake of contradiction, that there exist two valid choices for which block of $S$ contains $j$. That is, suppose that there are two possible extensions of the set partition corresponding to the pair $(\pi, P)$ that have $j$ in different blocks. By \Cref{lem:recoveralpha}, the content partition of $\alpha$ along with $\pi$ uniquely determines $\alpha$. In particular, as presented in the proof of the lemma, if $J$ is the block of $S$ containing $j$, then car $j$'s preference is $\pi^{-1}(k)$, where $k$ is the highest index in $J$. Thus, for there to be two possible extensions placing $j$ in different blocks, there must be two distinct possible preferences for car $j$, which we refer to as spots $A$ and $B$. We aim to show that since $A\neq B$, the final parking spot
		of car $j$ will be distinct in each of these scenarios. Ultimately, at most one of these scenarios corresponds to the outcome $\pi$.

		For the sake of readability, we summarize here the objects that will be used in the arguments that follow:
		\begin{itemize}
			\item $A$ and $B$ : the two distinct possible preferred spots for car $j$.

			\item $G_{m}$ : the set of distinct preferred spots among the first $m$ cars.

			\item $g_{m}:=|G_{m}|$ : the number of distinct preferred spots among the first $m$ cars.

			\item $h_{m}:=m-g_{m}$ : in the partial parking outcome after step $m$, this is the number of cars parked in thus-far undesired spots.

			\item $f(m)$ : the $h_{m}$-th lowest-numbered spot in $[n]\setminus G_{m}$.
		\end{itemize}
		Since neither $A$ nor $B$ could \textit{a priori} be ruled out as car $j$'s preferred spot, there must already be a car indexed by an element $a>j$ of $[n]$ corresponding to $))$ whose preference is known to be spot $A$, and a car indexed by an element $b>j$ of $[n]$ corresponding to $))$ whose preference is spot $B$. Furthermore, there is no car indexed by an element of $[n]$ greater than $j$ and corresponding to $(($ whose preference is spot $A$ or spot $B$.

		Regardless of whether car~$j$'s true preference is spot $A$ or spot $B$ (or something else), we know that in the partial parking outcome (see \Cref{def:partial parking outcome}) after step $j$, both spots $A$ and $B$ are occupied. Namely, in the case that car $j$ does not prefer spot $A$ (respectively, spot $B$), there exists an $\ell\in[j-1]$ corresponding to $(($ so that car $\ell$ prefers spot $A$ (respectively, spot $B$). Regardless of the preferences for the cars in $[j-1]$, the following holds: If car~$j$ initially parks in spot $A$, it will never occupy spot $B$ even after bumping, as that spot is already occupied. Similarly, if it initially parks in spot $B$, it will never occupy spot $A$.

		We consider the scenarios where car $j$ prefers either spot $A$ or spot $B$ and compare how car $j$ is affected by the parking of subsequent cars. In particular, we will show that after any possible sequence of bumps corresponding to the preferences of the cars with indices greater than $j$, there is no partial parking outcome where car $j$ is parked in the same spot if it originally preferred spot $A$ vs. spot $B$. Thus, at most one of these scenarios corresponds to our desired outcome $\pi.$

		For $m>j$, recall $G_{m}$ is the set of distinct preferred spots among the first $m$ cars and $g_{m}=|G_{m}|$. Note that $G_{m}$ and consequently $g_{m}$ do not depend on car $j$'s original preference, nor on any of the thus-far undetermined preferences of lower-indexed cars. In particular, as we already know the preferences of every car with an index corresponding to $))$ or $()$, we know the entire set of spots preferred by at least one car. As we also know the preferences of every car with an index greater than $j$, and hence of every car with an index greater than $m$, we also know which spots are desired by a car with an index greater than $m$ corresponding to $(($ or $()$ and thus by no car of lower index. If car $m$ performs a bump, we know that there is a repeated
		preference among the first $m$ cars so $g_{m}<m$ and the bumped car is bumped into the $(m-g_{m})$-th lowest-numbered spot in $[n]\setminus G_{m}$. Recall that $h_{m}:=m-g_{m}$ and $f(m)$ is the $h_{m}$-th lowest-numbered spot in $[n]\setminus G_{m}$, so any car bumped by car $m$ is bumped into spot $f(m)$.

		Suppose that the first time car~$j$ occupies the same spot in both scenarios is after car $p'$ parks and suppose without loss of generality that car $p'$ bumped car $j$ in the scenario where car $j$ preferred spot $B$. Let $p<p'$ be the index of the car that performed the most recent bump of car~$j$ in the scenario where car $j$ preferred spot $A$. Any car bumped by car $p$ is bumped to spot $f(p)$ while any car bumped by car $p'$ is bumped to spot $f(p')$ so $f(p)=f(p')$.

		Thus, the $h_{p}$-th lowest-numbered spot among those not desired by any of the first $p$ cars, (i.e., in $[n]\setminus G_{p}$), is the same as the $h_{p'}$-th lowest-numbered spot among those not desired by any of the first $p'$ cars. Note that for a given $t$ and $k_{1}<k_{2}$, the $t$-th lowest-numbered spot in $[n]\setminus{G_{k_2}}$, if such a spot exists, has an index at least as large as the $t$-th lowest-numbered spot in $[n]\setminus{G_{k_1}}$. By considering $t=h_{p}, k_{1}=p, k_{2}=p'$, we conclude that $h_{p'}\le h_{p}$.

		Furthermore, since $h_{i}$ is nondecreasing, we have in general that $f(i)$ is nondecreasing, a fact we use shortly. Also, since $p'>p$, we have $g_{p'}-g_{p}\le p'-p$, with equality only obtained if cars $p+1,\dots,p'$ all park in previously undesired spots. This rearranges to $p-g_{p}\le p'-g_{p}'$ or equivalently $h_{p}\le h_{p'}$. Thus $h_{p}=h_{p}'$, meaning that cars $p+1,\dots,p'$ all park in previously undesired spots, and in particular, none of them prefers a spot desired by any of the first $p$ cars. However, in either scenario, note that the first car to bump car~$j$ has the same preference as car $j$. Since $j<p$, it means the preference of the first car to bump car $j$ is among the preferences of the first $p$ cars, so in the scenario where car $j$ preferred spot $B$, car
		$p'$ cannot be the first car to bump car~$j$.

		Suppose car $q<p'$ is the car that most recently bumped car~$j$ before car $p'$ did, in the scenario where car $j$ preferred spot $B$. In particular, car $q$ bumps car $j$ to spot $f(q)$. Thus, car $p'$'s preferred spot is $f(q)$. As every car from car $p+1$ through car $p'$ parked in a spot undesired by the first $p$ cars, it means none of the first $p$ cars preferred spot $f(q)$.

		Let $l$ denote the common value $h_{p'}=h_{p}$. Thus, $f(p)$ is the $l$-th lowest-numbered spot in $[n]\setminus G_{p}$ and $f(p')$ is the $l$-th lowest-numbered spot in $[n]\setminus G_{p'}$. As these are equivalent, none of the $l$ lowest-numbered spots in $[n]\setminus G_{p}$ are desired by any of the cars $p+1,\dots, p'$. As $f(q)$ is desired by one of the cars $p+1,\dots, p'$, namely car $p'$, it means $f(q)$ is not among the $l$ lowest-numbered spots in $[n]\setminus G_{p}$. However, $f(q)$ is undesired by the first $p$ cars so there are at least $l$ lower-numbered spots also undesired by the first $p$ cars. In particular, spot $f(p)$ is the $l$-th lowest-numbered spot undesired by the first $p$ cars, so $f(q)>f(p)=f(p')$. Because $f$ is nondecreasing, this implies $q>p'$, a
		contradiction. Thus, for every pair $(\pi, P)$, with $\pi$ a permutation and $P$ a Dyck path, there is at most one set partition $S$ mapping to $P$, such that $S=\calS(\alpha)$ for some $\alpha$ with $\Out_{B}(\alpha)=\pi$.
	\end{proof}
	We are now ready to establish the proof of \Cref{prop:justupper}.

	\begin{proof}
		[Proof of \Cref{prop:justupper}] Consider the map from $\Out_{B}^{-1}(\pi)$ to the set of Dyck paths of semilength $n$ defined by sending $\alpha\in\Out_{B}^{-1}(\pi)$ to the pathification of $S(\alpha)$ (see Definition~\ref{def:dyckpath}). By \Cref{lem:onePrefListperP}, this map is injective. Thus $|\Out_{B}^{-1}(\pi)|$ is bounded above by the number of Dyck paths of semilength $n$, which is $C_{n}$.\qedhere
	\end{proof}

	Next we consider the case where the inequality in \Cref{prop:justupper} is strict.
	\begin{proposition}
		\label{thm:strictupper} Let $L=\{12\cdots (n-1)n,\  12\cdots(n-2)n(n-1), \ 12\cdots(n-3)n(n-2)(n-1)\}\subseteq \Sym_{n}$. Then for all $\pi\in\Sym_{n}\setminus L$, we have:
		\[
			\left|\Out_{B}^{-1}\left(\pi\right)\right|<C_{n}.
		\]
	\end{proposition}

	\begin{proof}
		For a given permutation $\pi\in\Sym_{n}$ and a Dyck path $P$ of semilength $n$, we established in \Cref{lem:onePrefListperP} that there is at most one preference list $\alpha$ with outcome $\pi$ whose content partition $\calS(\alpha)$ has pathification $P$. We show that for all $\pi \in \Sym_{n}\setminus L$, there is at least one Dyck path where there is no such preference list. To this end, consider the following Dyck paths, where a ``pairing'' refers to steps $2i-1$ and steps $2i$ for $1\le i\le n$:

		\begin{itemize}
			\item $Q=(()()(\dots)())$ i.e., the Dyck path that contains just one pairing $(($ at the beginning, just one pairing $))$ at the end, and no $()$ pairing, and

			\item $R=(()()(\dots)())()$ i.e., the Dyck path that contains just one pairing $(($, one pairing $))$, and one pairing $()$, where $(($ appears at the beginning, and the pairings $))$ and $()$ appear in this order at the end.
		\end{itemize}

		We will show that for every preference list $\alpha$ corresponding to $Q$, $\Out(\alpha)$ must be of the form $12 \cdots(r-1)n(r)(r+1)\cdots(n-1)$ for some $r\in[n]$. This shows that for $\pi$ not of this form, the path $Q$ suffices.

		For the remaining $\pi,$ we will show that if $\pi=\Out(\alpha)$ for some $\alpha$ whose content partition corresponds to $R$, then $\Out(\alpha)=12 \cdots(r-1)n(r)(r+1)\cdots(n-1)$ must have $r\geq n-2$, thus narrowing our choice of $\pi$ down to $L=\{12\cdots (n-1)n,\  12\cdots(n-2)n(n-1), \ 12\cdots(n-3)n(n-2)(n-1)\}$.

		Consider the Dyck path $Q=(()()(\dots)())$. Then, $Q=P_{S}$ if and only if $S=\{[n]\},$ the trivial set partition. If the content partition of some preference list corresponds to this Dyck path, it is a preference list where every car has the same preference. If every car prefers spot $r$, the resulting permutation has car $n$ in spot $r$ and the other cars in increasing order. Thus, only permutations of the form $12 \cdots(r-1)n(r)(r+1)\cdots(n-1)$ for some $r\in[n]$ could possibly have a preference list whose content partition corresponds to $Q$. This implies that there are strictly fewer than $C_{n}$ preference lists resulting in any other outcome.

		Next, consider the Dyck path $R=(()()(\dots)())()$. Then, $R=P_{S}$ if and only if $S$ is the set partition $\{\{1,2,\dots,n-1\},\{n\}\}$. If a preference list has this content partition and has outcome $12 \cdots(r-1)n(r)(r+1)\cdots(n-1)$, then car $n$ prefers spot $r$ and no other car prefers spot $r$. Hence, car $n$ must not have bumped car $n-1$, as they do not share a preference.

		Thus, car $n-1$ must have preferred the spot it ends up parked in which is spot $n$ when $r\le n-1$. Consequently, for $r\le n-1$, every car but car $n$ prefers spot $n$. However, if $r\le n-2$, this means that in the partial parking outcome after step $n-1$, car $r$ is in spot $r$ and spot $n-1$ is empty. Thus, when car $n$ parks, car $r$ is bumped to spot $n-1$. For this to actually result in the outcome $12 \cdots(r-1)n(r)(r+1)\cdots(n-1)$, we must have $r=n-2$. Thus, aside from possibly the permutations $12\cdots n, 12\cdots(n-2)(n)(n-1),$ and $12\cdots(n-3)n(n-2)(n-1)$, there are no others where two preference lists resulting in that outcome have content partitions corresponding to $Q=(()()(\dots)())$ and $R=(()()(\dots)())()$. Thus, every permutation outside $L$ is obtained as
		an outcome by strictly fewer than $C_{n}$ preference lists.
	\end{proof}

	Next, we make the following definition.
	\begin{definition}
		Suppose car $j$ is parked in spot $i$ on the street. If car $j$ is bumped out of spot $i$ and reparks in a spot numbered $k>i$, we say that car $j$ was \textit{bumped forward}. If instead car $j$ is bumped out of spot $i$ and reparks in a spot numbered $k<i$, we say that car $j$ was \textit{bumped backward}.
	\end{definition}

	\begin{remark}
		\label{rmk:bumpbackonce} Note that under the backgammon parking scheme, whenever a car gets bumped, it gets bumped to a spot that has no open spots preceding it. Consequently, if the car is bumped again, it must be bumped forward. Therefore, aside from possibly the first time it is bumped, a car is never bumped backward, so no car is ever bumped backward more than once.
	\end{remark}

	\begin{proposition}
		\label{lem:combined} Let $L=\{12\cdots (n-1)n,\  12\cdots(n-2)n(n-1), \ 12\cdots(n-3)n(n-2)(n-1)\}\subseteq \Sym_{n}.$ Then, for all $\pi\in L$, we have $|\Out_{B}^{-1}(\pi)|=C_{n}$.
	\end{proposition}

	\begin{proof}
		Fix $\pi\in L$. Recall that a set partition is said to be non-nesting if we cannot find four elements $i < j < j' < i'$ with $i,i'$ in one block and $j,j'$ in another block such that there is no other $k$ in the block with $i,i'$ satisfying $i<k<i'$. We will show that for each non-nesting set partition $S$ of $[n]$, there exists a preference list $\alpha_{S}$ with content partition $\mathcal{S}(\alpha_{S})=S$ and outcome $\Out(\alpha_{S})=\pi$. Consequently, the number of preference lists yielding $\pi$ is at least the number of non-nesting partitions of $[n]$, which is $C_{n}$ (see \seqnum{A000108}). Along with the upper bound in \Cref{prop:justupper}, this establishes $|\Out_{B}^{-1}(\pi)|=C_{n}$.

		For a given non-nesting set partition $S$, we construct our desired $\alpha_{S}\coloneqq(a_{1},a_{2},\ldots,a_{n})$ as follows. Let $J$ be a block of $S$, and let $k=\max(J)$ be the highest index in $J$. For each index $j$ in $J$, assign the preference $a_{j}=\pi^{-1}(k)$. `` Any car whose index is highest in its block is the last car with that preference to arrive, and therefore parks in its preferred spot. Thus any such car $k$, including car $k=n$, parks in the spot $\pi^{-1}(k)$ as desired. We must now check that the remaining cars park according to $\pi$ and it suffices to check that they park in increasing order.

		First, we establish that for the cars that are bumped, their respective first bumps occur in increasing order. Suppose otherwise and that there are indices $i<j$ and $j'<i'$ such that car $i$ is bumped for the first time by car $i'$ and car $j$ is bumped for the first time by car $j'$. As $j'$ must be greater than $j$, this yields $i<j<j'<i'$. Furthermore, cars $i,i'$ must have the same preference and no car $k$ with $i<k<i'$ can share that preference, as otherwise car $k$ would have bumped car $i$ earlier than car $i'$. Similarly, cars $j,j'$ must have the same preference and no car $k$ with $j<k<j'$ can share that preference. Thus, the parts of $S$ containing $\{i,i'\}$ and $\{j,j'\}$ would violate the fact that $S$ is non-nesting.

		As the first bumps occur in increasing order, the only way the bumped cars will not park in increasing order is if some car is bumped more than once. We show that this is only possible when $\pi=12\cdots(n-3)(n)(n-2)(n-1)$ and in that case, at most one car (in particular, the car indexed $n-2$) can be bumped twice.

		For $\pi\in\{12\cdots n,12\cdots(n-2)n(n-1)\}\subset L$, we show that no car is bumped more than once. Let $k=\max(J)$ for some block $J$ of $S$. The only time car $k$ prefers a spot whose index is smaller than $k$ is when $\pi=12\cdots(n-2)n(n-1)$ and $k=n$, in which case car $k=n$ prefers spot $n-1$. Lower-numbered indices in the same block correspond to cars with the same preference, so any car $j$'s preferred spot is at least $j$ except when $\pi=12\cdots(n-2)n(n-1)$ and $j=n$. If a car $i$ is bumped for the first time by car $j$, it necessarily is bumped into one of the first $j$ spots. Then, in the case of $12\cdots n$, every subsequent car has a preferred spot of index at least $j+1$ so car $i$ cannot possibly be bumped again. The only possible exception for $\pi=12\cdots(n-2
		)n(n-1)$ is if car $i$ is first bumped into spot $n-1$ and then bumped a second time by car $n$. However, the only way a car could be bumped into spot $n-1$ by one of the first $n-1$ cars is if it is bumped by car $n-1$ and none of the first $n-1$ cars prefer spot $n$. Based on how $\alpha_{S}$ is constructed, car $n-1$ would prefer spot $n$ if $n-1$ is the highest index in its block of $S$. Thus, it suffices to consider when $n-1$ is not the highest index in its block of $S$, which would mean cars $n-1$ and $n$ have the same preference in $\alpha_{S}$. This would mean car $n$ bumps car $n-1$ and consequently cannot bump the car already bumped by car $n-1$. Thus, no car is bumped twice, as desired.

		Now, we fix $\pi=12\cdots(n-3)(n)(n-2)(n-1)$ and show that only one car can possibly be bumped twice. Let $k=\max(J)$ for some block $J$ of $S$. Then, the only time car $k$ has a preference in $\alpha_{S}$ that is less than $k$ is when $k=n.$ In that case, car $k$ prefers spot $n-2$, so for any car $j\le n-2$, the index of car $j$'s preferred spot is at least $j$. Cars $n$ and $n-1$ each have a preferred spot with index at least $n-2$. If a car $i$ is bumped for the first time by a car $j\le n-3$, it is necessarily bumped into one of the first $j$ spots. Then, every subsequent car has a preferred spot of index at least $j+1$ so car $i$ cannot possibly be bumped again. Thus, for a car to be bumped multiple times, it must be bumped by at least two of the three highest-numbered cars. Therefore,
		at most one car, say car $c$, is bumped at least twice.

		The only way the bumped cars do not end up parked in increasing order is if some car $d$ is bumped after the first time car $c$ is bumped but before the last time car $c$ is bumped. The first time car $c$ is bumped is by a car with index at least $n-2$, so this is only possible if car $c$ is bumped by cars $n-2$ and $n$, and car $d$ is bumped by car $n-1$. For this to happen, car $n-2$ must bump car $c$ into spot $n-2$, meaning none of the first $n-2$ cars have a preference outside the first $n-2$ spots. Also, since car $n$ does not bump car $n-1$, it means that $n-1$ is the highest index in its block of $S$, so car $n-1$ must prefer spot $n$. This gives a contradiction, since none of the first $n-2$ cars prefer spot $n$ but car $n-1$ is supposed to park in a spot preferred by car
		$d$. Thus, even if some car ends up being bumped twice, the bumped cars among $1,\dots,n-1$ still end up in increasing order. Therefore, for each non-nesting set partition $S$, the preference list $\alpha_{S}$ yields the outcome $12\cdots(n-3)n(n-2)(n-1)$, as desired. \qedhere
	\end{proof}

	\Cref{thm:at most Catalan} follows directly from \Cref{prop:justupper}, \Cref{thm:strictupper}, and \Cref{lem:combined}. We are now ready to show that for any permutation $\pi\in\Sym_{n}$, there are at least $2^{n-1}$ preference lists with backgammon outcome equal to $\pi$.

	\BackgammonLowerBound*

	\begin{proof}
		Define $L_{n}:=\displaystyle\min_{\pi\in\Sym_n}|\Out_{B}^{-1}(\pi)|$ to be the smallest number of preference lists that can lead to an outcome $\pi\in\Sym_{n}$. We will show that $L_{n}\geq 2^{n-1}$. By direct computation, we determine the following values: $L_{1}= 1$, $L_{2}= 2$, $L_{3}= 4$, $L_{4}= 8,$ and $L_{5}=17$.

		For $n\geq 6$, fix an outcome $\pi\in S_{n}$. For each $1\le i \le n$, there exists a partial parking outcome $\beta_{i}$ consisting of $n-1$ cars occupying all spots except the spot $i$, such that parking car $n$ in spot $\pi^{-1}(n)$ will yield outcome $\pi$. For each such partial parking outcome $\beta_{i}$, we consider preference lists satisfying all the following properties:
		\begin{itemize}
			\item car $n$ prefers spot $\pi^{-1}(n)$,

			\item the partial parking outcome after step $n-1$ is $\beta_{i}$,

			\item any car $j<n$ occupying a spot greater than $i$ in the partial parking outcome after step $n-1$ prefers that spot, and

			\item any car $j<n$ occupying a spot less than $i$ in the partial parking outcome after step $n-1$ prefers a spot among $\{1,2,\dots,i-1\}$.
		\end{itemize}
		As all such preference lists have outcome $\pi$, counting them will yield a lower bound for $|\Out_{B}^{-1}(\pi)|$.

		By construction, for each such preference list, the only bumps before car $n$ parks occur between the $i-1$ cars that are in spots less than $i$ in the partial parking outcome after step $n-1$. Thus, for $i\ge 2$, the number of ways for this to happen and to produce the outcome $\pi$ is at least the minimum number of preference lists for $i-1$ cars that yield a given outcome of length $i-1$. Also, including the case where $i=1$, which occurs in precisely one way, we get a lower bound
		\[
			L_{n}\ge 1+\sum_{i=2}^{n}L_{i-1}.
		\]
		By induction, this gives $L_{n}\ge 2^{n-1}$. Additionally, $L_{5}=17$, so for $n\ge 5$, we have $L_{n}>2^{n-1}$.
	\end{proof}

	\begin{remark}
		Since we have $n^{n}$ possible preference lists divided among $n!>(n/e)^{n}$ outcomes, some outcome must arise from fewer than $e^{n}$ preference lists. This means that the smallest number of preference lists giving a particular outcome of length $n$ is strictly less than $e^{n}$.
	\end{remark}

	\subsection{Decreasing Permutations}
	\label{sec:dec perms in backgammon} There are several other permutations $\pi$ for which we can compute the exact size of $\Out_{B}^{-1}(\pi)$. These results are related to Motzkin numbers, which first appeared in the work of Motzkin \cite{motzkin}, and have since received considerable attention in combinatorics \cite{DS77}. For our purposes, we recall the definition of Motzkin paths, which are enumerated by the Motzkin numbers.
	\begin{definition}
		\label{def:Motzkin paths} Let $n \in \NN$. A \textit{Motzkin path of length $n$} is a lattice path consisting solely of horizontal steps $(1,0)$ (denoted $\tH$), upward diagonal steps $(1,1)$ (denoted $\U$), and downward diagonal steps $(1,-1)$ (denoted $\D$), which begins at $(0,0)$ and ends at $(n,0)$, and which does not fall below the $x$-axis. We let $\mathcal{M}_{n}$ denote the set of Motzkin paths of length $n$ and we let $M_{n}=|\mathcal{M}_{n}|$.
	\end{definition}
	We recall that $M_{n}$ is called the $n$-th Motzkin number and satisfies the formula
	\begin{align}
		\label{eq:motzkin}M_{n}=\sum_{k=0}^{\lfloor\frac{n}{2}\rfloor}\binom{n}{2k}C_{k},
	\end{align}
	where $C_{k}$ is the $k$-th Catalan number. The sequence $\{M_{n}\}_{n\geq 1}$ is OEIS sequence \seqnum{A001006}. \Cref{fig:motzkin} illustrates the nine Motzkin paths of length $4$.

	\begin{figure}
		[h!] \centering \resizebox{\textwidth}{!}{ \begin{tikzcd}\draw[step=1cm,gray,very thin] (0,0) grid (4,2); \draw[ultra thick, red] (0,0)--(4,0); \node at(2,-1){\Huge{\tH\tH\tH\tH}};\end{tikzcd} \qquad \begin{tikzcd}\draw[step=1cm,gray,very thin] (0,0) grid (4,2); \draw[ultra thick, red] (0,0)--(2,0)--(3,1)--(4,0); \node at(2,-1){\Huge{\tH\tH\U\D}};\end{tikzcd} \qquad \begin{tikzcd}\draw[step=1cm,gray,very thin] (0,0) grid (4,2); \draw[ultra thick, red] (0,0)--(1,0)--(2,1)--(3,0)--(4,0); \node at(2,-1){\Huge{\tH\U\D\tH}};\end{tikzcd} \qquad \begin{tikzcd}\draw[step=1cm,gray,very thin] (0,0) grid (4,2); \draw[ultra thick, red] (0,0)--(1,1)--(2,0)--(3,0)--(4,0); \node at(2,-1){\Huge{\U\D\tH\tH}};\end{tikzcd} \qquad \begin{tikzcd}\draw[step=1cm,gray,very thin] (0,0) grid (4,2); \draw[ultra thick, red] (0,0)--(1,0); \draw[ultra thick, red] (1,0)--(2,1)--(3,1)--(4,0); \node at(2,-1){\Huge{\tH\U\tH\D}};\end{tikzcd} \qquad \begin{tikzcd}\draw[step=1cm,gray,very thin] (0,0) grid (4,2); \draw[ultra thick, red] (0,0)--(1,1)--(2,1)--(3,0)--(4,0); \node at(2,-1){\Huge{\U\tH\D\tH}};\end{tikzcd} \qquad \begin{tikzcd}\draw[step=1cm,gray,very thin] (0,0) grid (4,2); \draw[ultra thick, red] (0,0)--(1,1)--(2,1)--(3,1)--(4,0); \node at(2,-1){\Huge{\U\tH\tH\D}};\end{tikzcd} \qquad \begin{tikzcd}\draw[step=1cm,gray,very thin] (0,0) grid (4,2); \draw[ultra thick, red] (0,0)--(1,1)--(2,0)--(3,1)--(4,0); \node at(2,-1){\Huge{\U\D\U\D}};\end{tikzcd} \qquad \begin{tikzcd}\draw[step=1cm,gray,very thin] (0,0) grid (4,2); \draw[ultra thick, red] (0,0)--(2,2)--(4,0); \node at(2,-1){\Huge{\U\U\D\D}};\end{tikzcd} }
		\caption{Illustration
		of
		the
		nine
		distinct
		Motzkin
		paths
		of
		length
		$4$.}
		\label{fig:motzkin}
	\end{figure}
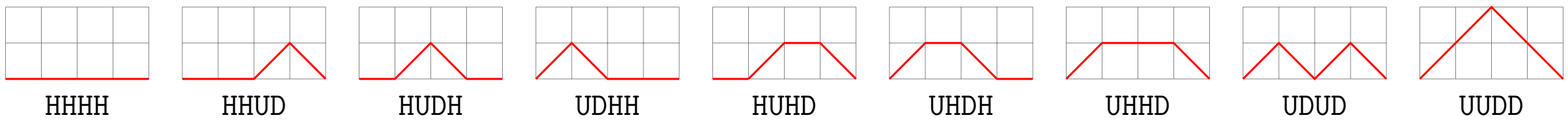

	We recall the following definition of permutation patterns:
	\begin{definition}
		Fix $m\leq n$ and let $\pi\in \Sym_{n}$ and $\sigma\in \Sym_{m}.$ Then, $\pi$ is said to contain a $\sigma$-\textit{pattern}, if some subsequence of the entries of $\pi$ has the same relative order as all of the entries of $\sigma$.
	\end{definition}
	For instance, permutation $\pi$ contains a $213$-pattern whenever $\pi$ has three entries $x,$ $y$, and $z$ that appear within $\pi$ in the order $x\cdots y\cdots z$ but whose values are ordered as $y < x < z$, the same as the ordering of the values in the permutation $213$.

	\begin{example}
		For example, $\pi=536421$ contains several instances of a $231$-pattern, including the subsequences $564$ and $341$.
	\end{example}
	
    We now state the main result of this section. \thmMotzkin* The proof of \Cref{thm:motzkin} relies on the following technical lemmas.
	
    \begin{lemma}
		\label{lem: nobumptwice} If the backgammon outcome is $n(n-1)\cdots 21$, then no car is bumped forward twice. Consequently, no car is bumped more than twice.
	\end{lemma}
	\begin{proof}
		Suppose for the sake of contradiction that car $i$ is bumped forward by cars $j$ and $k$ with $i<j<k$, and consider the partial parking outcome after step $k$. Since car $i$ was bumped forward by car $k$, immediately after that bump, car $k$ is in a lower-numbered spot than car $i$. By Remark~\ref{rmk:bumpbackonce}, car $i$ could not have been bumped backward by any cars between $j$ and $k$. Thus, if car~$j$ has not been bumped, then car $j$ is in a lower-numbered spot than car $k$. If car~$j$ has been bumped, then the highest numbered car between $j$ and $k$ which shares $j$'s preference is in a lower-numbered spot than car $k$. Either way, the partial parking outcome after step $k$ contains a $231$-pattern with no empty spots before the $1$. We now show that once we have achieved
		a $231$-pattern with no empty spots before the $1$, after each subsequent step, our partial parking outcome contains at least one of the following:
		\begin{itemize}
			\item a $231$-pattern with no empty spots before the $1$, or

			\item a $312$-pattern with no empty spots before the $2$.
		\end{itemize}

		This would imply that we can never achieve the outcome where the cars are parked in decreasing order. If we have a $231$-pattern with no empty spots before $1$, we continue to have this unless one of the three cars representing $2,3,$ or $1$ in the pattern is bumped by a higher car representing $4$. This results in having a partial parking outcome containing one of the patterns $4312, 2413$, or $2341$, where again there are no empty spots before the last of these four cars. Both patterns $2413$ and $2341$ contain a $231$-pattern, while $4312$ contains a $312$-pattern. If we have a $312$-pattern with no empty spots before the $2$, we continue to have this unless one of the three cars representing $3,1,$ or $2$ is bumped by a higher car representing $4$. This results in $4123, 3421$,
		or $3142$ with again no empty spots before the last of these four cars. Both $4123$ and $3142$ contain a $312$-pattern, while $3421$ contains a $231$-pattern. Thus, when the final backgammon outcome is $n(n-1)\cdots 21$, no car is bumped forward more than once. Note that in general, by Remark~\ref{rmk:bumpbackonce}, no car is ever bumped backward more than once. Thus, if the final outcome is $n(n-1)\cdots 21$, no car is bumped more than twice.
	\end{proof}
	For the sake of completion, we prove two recurrences for the triangle of Motzkin ballot numbers. The recurrence in \Cref{rec:OEIS} is stated in \cite[\seqnum{A091836}]{OEIS} without proof. These recurrences are used in our proof of \Cref{thm:motzkin}.
	\begin{lemma}
		\label{lem:OEIS recurrence} Let $t(n,k)$ denote the number of Motzkin paths with $n$ steps that touch $k-1$ points on the $x$-axis besides the beginning and the end. For $n\ge 2$ and $k\geq 2$, $t(n,k)$ satisfies the following recurrence:
		\begin{align}
			t(n,k) & =t(n,k-1)+t(n-1,k-1)-t(n-1,k-2)-t(n-2,k-2). \label{rec:OEIS}
		\end{align}
		Furthermore, when $n\ge 3$ and $k\ge 1$, we have:
		\begin{align}
			t(n,k) & =t(n-1,k-1)+\sum_{i=k-1}^{n-2}t(n-2,i).\label{rec:second recurrence for Motzkin}
		\end{align}
	\end{lemma}
	\begin{proof}
		Note that $t(n,k)$ is understood to be $0$ if $k>n$ or if $k=0$ and $n>0$. For induction purposes, we define $t(0,0)\coloneqq 1$. Let $g(x)=\sum_{i=0}^{\infty}M_{i-1}x^{i}$ with the convention that $M_{-1}\coloneqq 1$. For $k\ge 1$, a Motzkin path which touches the $x$-axis at $k-1$ points besides the beginning and the end can be viewed as the composition of $k$ Motzkin paths which each only touch the $x$-axis at their respective endpoints. A Motzkin path of length $i\ge2$ with no intermediate points on the $x$-axis starts with a $\U$ step and ends with a $\D$ step. Having no other intermediate points on the $x$-axis means that never in the middle $i-2$ steps does the path go below where it starts and finishes those $i-2$ steps. Thus, there are $M_{i-2}$ Motzkin paths of length $i$
		which hit the $x$-axis only at their endpoints. (Note that the previous statement also holds for $i=1,$ since $M_{-1}=1$). Therefore, for $(n,k)\ne(0,0)$, $t(n,k)$ is the sum of the products $M_{i_1-2}M_{i_2-2}\cdots M_{i_k-2}$ over all tuples $(i_{1},i_{2},\dots,i_{k})$ of positive integers satisfying $i_{1}+i_{2}+\dots+i_{k}=n$. As $M_{i-2}$ is the coefficient of $x^{i-1}$ for $g(x)$, we obtain that $t(n,k)$ is the $[x^{n-k}]$ coefficient of $g(x)^{k}$ for $k\ge 1$. Note that this still holds for $k=0,n\ge 0$ since for $n\ge 1$, $[x^{n}]g(x)^{0}=0=t(n,0)$ and for $n=0$, $[x^{n}]g(x)^{0}=1=t(0,0)$.

		After reindexing, the recurrence we wish to show for $n,k\ge 0$ is
		\begin{align}
			\label{eq:four term rec for t}t(n+2,k+2)=t(n+2,k+1)+t(n+1,k+1)-t(n+1,k)-t(n,k).
		\end{align}

		If $k>n+1$, then \cref{eq:four term rec for t} simply becomes $0=0+0-0-0$. If $k=n+1$, then \cref{eq:four term rec for t} becomes $0=t(n+2,n+2)+0-t(n+1,n+1)-0$. This is true since $t(m,m)=1$ for all $m\ge 1$ because the only Motzkin path of length $m$ touching the $x$-axis at $m-1$ points besides the beginning and the end is the path consisting solely of $\tH$ steps. Thus we may assume $n\ge k$. In this case, we wish to show that the coefficient $[x^{n-k}]$ of $g(x)^{k+2}$ is equal to the sum of the coefficients $[x^{n-k}]$ and $[x^{n-k+1}]$ for $(g(x)-1)g(x)^{k}$.

		Let $j=n-k$. Since $n\ge k$, we have $j\ge 0$, and $n=j+k$. Rewriting the desired identity in terms of $j$, we obtain
		\[
			[x^{j}]g(x)^{k+2}=[x^{j}]\left((g(x)-1)g(x)^{k}\right)+[x^{j+1}]\left((g(x)-1)g(x)^{k}\right).
		\]
		Thus, it suffices to prove that for all $k\ge 0$, the above statement holds for every $j\ge 0$. We proceed by induction on $k$. For the base case $k=0$, we split into $j=0$ and $j\ge 1$. To handle $j=k=0$, we note that
		\[
			[x^{0}]g(x)^{2}=1=0+1=[x^{0}](g(x)-1)+[x^{1}](g(x)-1).
		\]
		For $j\ge 1$, we must show that $[x^{j}]g(x)^{2}= 2M_{j-1}+ \sum_{i=1}^{j-1}M_{i-1}M_{j-i-1}$ is equal to $[x^{j}](g(x)-1)+[x^{j+1}](g(x)-1)= M_{j}+ M_{j-1}$. This follows from the well-known identity \cite[page 90]{BERNHART199973}:
		\[
			M_{n}=M_{n-1}+\sum_{i=1}^{n-1}M_{i-1}M_{n-i-1}.
		\]

		The induction hypothesis is that for some $k$, we have for all $j\ge 0$, that
		\[
			[x^{j}]g(x)^{k+2}=[x^{j}]\left((g(x)-1)g(x)^{k}\right)+[x^{j+1}]\left((g(x)-1)g(x)^{k}\right)
		\]
		Inducting on $k$, we obtain for $j\ge 0$,
		\begin{align*}
			[x^{j}]g(x)^{k+3} & =\sum_{i=1}^{j+1}([x^{i-1}]g(x))([x^{j+1-i}]g(x)^{k+2})                                   \\
			                  & =\sum_{i=1}^{j+1}([x^{i-1}]g(x)) \left(([x^{j+1-i}]+[x^{j+2-i}])((g(x)-1)g(x)^{k})\right) \\
			                  & =\sum_{i=1}^{j+1}([x^{i-1}]g(x)) ([x^{j+1-i}]((g(x)-1)g(x)^{k}))                          \\
			                  & \qquad\qquad+ \sum_{i=1}^{j+1}([x^{i-1}]g(x)) ([x^{j+2-i}]((g(x)-1)g(x)^{k}))             \\
			                  & =[x^{j}]((g(x)-1)g(x)^{k+1})                                                              \\
			                  & \qquad\qquad+ \sum_{i=1}^{j+2}([x^{i-1}]g(x)) ([x^{j+2-i}]((g(x)-1)g(x)^{k}))             \\
			                  & =[x^{j}]((g(x)-1)g(x)^{k+1}) +[x^{j+1}]((g(x)-1)g(x)^{k+1}),
		\end{align*}
		where the penultimate equality follows from the fact that $[x^{0}]((g(x)-1)g(x)^{k})=0$. This completes the proof of \eqref{rec:OEIS}.

		Next, we claim that for $n\ge 3, k\ge 1$:
		\[
			t(n,k)=t(n-1,k-1)+\sum_{i=k-1}^{n-2}t(n-2,i).
		\]

		Recalling that for $i\geq 1$, $t(i,1)=M_{i-2}$ and $t(i,0)=0$, we have for $n\ge 3$, that
		\[
			t(n,1)=M_{n-2}=t(n-1,0)+\displaystyle\sum_{i=0}^{n-2}t(n-2,i).
		\]

		We now induct on $k$. Beginning with \Cref{rec:OEIS}, for $n\ge 2$ and $k\ge 2$ we have that
		\begin{align*}
			t(n,k) & =t(n,k-1)+t(n-1,k-1)-t(n-1,k-2)-t(n-2,k-2)                                           \\
			       & =\left(t(n-1,k-2)+\sum_{i=k-2}^{n-2}t(n-2,i)\right)+t(n-1,k-1)-t(n-1,k-2)-t(n-2,k-2) \\
			       & =\sum_{i=k-2}^{n-2}t(n-2,i)+t(n-1,k-1)-t(n-2,k-2)                                    \\
			       & =t(n-1,k-1)+\sum_{i=k-1}^{n-2}t(n-2,i),
		\end{align*}
		as desired.
	\end{proof}

	\begin{proof}
		[Proof of \Cref{thm:motzkin}] As before, let $t(n,k)$ denote the number of Motzkin paths with $n$ steps, which touch $k-1$ points on the $x$-axis besides the beginning and the end. For $n,k\ge 1$, let $f(n,k)$ denote the number of preference lists yielding outcome $n(n-1)\cdots 321$, where the longest initial sequence of distinct preferences has length $k$. Note that $f(n,k)=0$ for $k>n$.

		Next, we show that $f(n,k)=t(n,k)$ for all $n\ge k\ge 1$. First note that $f(1,1)=f(2,1)=f(2,2)=1$ so the statement holds for $n<3$. Thus, the total number of preference lists yielding $n(n-1)\cdots321$ is $M_{n}$ for $n=1,2$. To prove the desired equalities for $n\ge 3$, we now induct on $n$ with the induction hypothesis that $f(m,k)=t(m,k)$ for all $1\le k\le m\le n-1$. To show that $f(n,k)=t(n,k)$ for all $1\le k\le n$, it suffices to show that the following hold for each fixed $n$:
		\begin{claim}
			$f(n,1)=M_{n-2}=t(n,1)$.
		\end{claim}
		\begin{claim}
			For $k\ge 2$, $f(n,k)$ obeys the same recurrence as $t(n,k)$, namely
			\[
				f(n,k)=f(n-1,k-1)+\sum_{i=k-1}^{n-2}f(n-2,i).
			\]
		\end{claim}

		\begin{proof}
			[Proof of Claim 1]

			Observe that $f(n,1)$ counts sequences yielding outcome $n(n-1)\cdots321$, where cars $1$ and $2$ have the same preference. If this preference is spot $1$, then car $2$ bumps car $1$ forward into spot $2$. Eventually, car $1$ will need to be bumped forward into spot $n$ but by \Cref{lem: nobumptwice}, it cannot be bumped forward twice. Thus, car $1$ does not prefer spot $1$ and is thus bumped into spot $1$ by car $2$. In order to end up in spot $n$, car $1$ must be the last car bumped. Thus, car $1$ must still be parked in spot $1$ in the partial parking outcome after step $n-1$, so that it is bumped by car $n$. Thus no car from those numbered $2$ to $n-1$ prefers spot $1$ and also no car from those numbered $2$ to $n-1$ can prefer spot $n$ since that must be left open for car $1$
			to be bumped into. Therefore, cars $2$ to $n-1$ park into spots $2$ to $n-1$ such that they end up in decreasing order. Since none of these cars can be bumped to spot $1$, this is equivalent to cars $1$ to $n-2$ parking in spots $1$ to $n-2$ in decreasing order. This is $\sum_{k=1}^{n-2}f(n-2,k)$ which by the induction hypothesis is equal to $\sum_{k=1}^{n-2}t(n-2,k)=M_{n-2}$, as desired.
		\end{proof}

		\begin{proof}
			[Proof of Claim 2] Now for $n\ge 3, k\ge 2$, consider the ways in which car $1$ could end up in spot~$n$. Either
			\begin{itemize}
				\item[i)] car $1$ prefers spot $n$,

				\item[ii)] car $1$ prefers spot $1$ and is bumped forward by car $n$, or

				\item[iii)] car $1$ prefers something else, is bumped to the start of the street and parks in spot $1$, and then is later bumped forward by car~$n$.
			\end{itemize}

			Since by \Cref{rmk:bumpbackonce}, a car cannot be bumped backward more than once, and by \Cref{lem: nobumptwice}, a car cannot be bumped forward more than once, these are the only possibilities.

			For case i), no other car can prefer spot $n$, so restricting to cars $2$ through $n$, the length of the longest initial sequence of distinct preferences is $k-1$. Cars $2$ through $n$ park in spots $1$ through $n-1$ in decreasing order, which is equivalent to cars $1$ to $n-1$ doing the same thing, so this case contributes $f(n-1,k-1)$ to the count.

			For case ii), we must have $k\le n-1$ as cars $1$ and $n$ share a preference for spot $1$. Cars $2$ through $n-1$ must have preferences among spots $2$ through $n-1$ so that spot $n$ is left open for car $1$ to be bumped into. No other car prefers spot $1$ besides cars $1$ and $n$, so when restricting to cars $2$ through $n-1$, the length of the longest initial sequence of distinct preferences is $k-1$. The number of such preference lists for cars $2$ through $n-1$ that park them in decreasing order in spots $2$ through $n-1$ is then $f(n-2,k-1)$. Car $n$ has no choice but to park in spot $1$, so this case contributes $f(n-2,k-1)$.

			For case iii), no car besides car $n$ can prefer spot $1$, and car $1$ must be the first car bumped. Since the first $k$ preferences are distinct and the longest such prefix has length $k$, car $k+1$ is the first car sharing a preference with an earlier car. Since car $1$ must be the first bumped car, that repeated preference must be car $1$'s preference. Thus car $k+1$ initially bumps car $1$. Since car $1$ must be bumped at least twice, we have that $k\le n-2$. Restricting to just cars $2$ through $n-1$, the initial sequence of distinct preferences has length at least $k$. When it is $i$, there are $f(n-2,i)$ ways to park cars $2$ through $n-1$ in decreasing order in spots $2$ through $n-1$ as none of those cars can ever occupy spot $1$ nor spot $n$. The preference of car $1$
			must then match the preference of car $k+1$ and the preference of car $n$ is spot $1$. Thus, for each $k\leq n-2$, this contributes a total of $\sum_{i=k}^{n-2}f(n-2,i)$ to the count.

			Tabulating the three cases for $k\le n-1$, we get
			\[
				f(n,k)=f(n-1,k-1)+f(n-2,k-1)+\sum_{i=k}^{n-2}f(n-2,i)=f(n-1,k-1)+\sum_{i=k-1}^{n-2}f(n-2,i),
			\]
			while for $k=n$, we get
			\[
				f(n,k)=f(n-1,k-1)=f(n-1,k-1)+\sum_{i=k-1}^{n-2}f(n-2,i),
			\]
			thus establishing the recurrence for $k\ge 2$ and completing the proof that $f(n,k)=t(n,k)$.
		\end{proof}
		Therefore, the total number of preference lists yielding outcome $n(n-1)\cdots321$ is $\sum_{i=1}^{n}f(n,i)=\sum_{i=1}^{n}t(n,i)=M_{n}$, as desired.
	\end{proof}

	Next we show that the number of backgammon parking functions with two special outcomes depends on the Motzkin numbers.

	\differenceofadjacentmotzkin*

	\begin{proof}
		Let $P_{n}:=\left|\Out_{B}^{-1}\left(1(n)(n-1)\cdots32\right)\right|$, with $P_{1}=1$ and $P_{2}=2$. Now consider $n\ge 3$. In order for the cars to park in decreasing order, there are $n$ possibilities for the partial parking outcome after step $n-1$. It is either $0(n-1)(n-2)\cdots 1$, $(n-1)0(n-2)(n-3)\cdots 1$ or $(k)(n-1)(n-2)\cdots(k+1)0(k-1)(k-2)\cdots1$ for some $1\le k\le n-2$ where car $k$ is the car that will be bumped by car $n$ and which will ultimately park in the remaining open spot, which is denoted by a $0$.

		Notably, all the cars after the empty spot have lower indices than the cars before the empty spot, so no car parked before the empty spot could have had a preference after the empty spot. The cars parked after the empty spot must have preferred the spots they ended up in. Thus, there is only one way to obtain each of the partial parking outcomes $0(n-1)(n-2)\cdots 1$ and $(n-1)0(n-2)(n-3)\cdots 1$. For $1\le k\le n-2$, the number of ways for the cars $k$ through $n-1$ to park before the empty spot in the order $(k)(n-1)(n-2)(n-3)\cdots(k+1)$ is the same as $P_{n-k}$, the number of ways for the first $n-k$ cars to park in the order $1(n-k)(n-k-1)\cdots2$. Thus, the number of ways for the cars to park in decreasing order, which is given in \Cref{thm:motzkin} as $M_{n}$, is equal to $2
		+\sum_{k=1}^{n-2}P_{n-k}$. Therefore,
		\[
			M_{n+1}-M_{n}=\left(2+\sum_{k=1}^{n-1}P_{n+1-k}\right)-\left(2+\sum_{k=1}^{n-2}P_{n-k}\right)=\sum_{j=2}^{n}P_{j}-\sum_{j=2}^{n-1}P_{j}=P_{n},
		\]
		as desired.
	\end{proof}

	\moremotzkin*

	\begin{proof}
		Let $Q_{n}:=\left|\Out_{B}^{-1}\left(12n(n-1)\cdots3\right)\right|$ with $Q_{1}=1$, $Q_{2}=2$, and $Q_{3}=5$. Now consider $n\ge 4$. In order for the cars to park in the order $1(n)(n-1)\cdots32$, the partial parking outcome after step $n-1$ is either $(1)0(n-1)(n-2)\cdots32$, $0(1)(n-1)(n-2)\cdots32$, $12(n-1)(n-2)\cdots 30$, $1(n-1)0(n-2)(n-3)\cdots32$, or $1(k)(n-1)(n-2)\cdots(k+2)(k+1)0(k-1)(k-2)\cdots32$, for some $3\le k\le n-2$, where car $k$ is the car that will be bumped by car $n$ and which will ultimately park in the remaining open spot, which is denoted by a $0$.

		The partial parking outcomes $(1)0(n-1)(n-2)\cdots32$ and $1(n-1)0(n-2)(n-3)\cdots32$ can each be attained in $n-1$ ways as any car with index greater than $1$ is parked in its preferred spot while car $1$ either preferred spot $1$ or was bumped by any of the cars from $2$ through $n-1$. The partial parking outcome $0(1)(n-1)\cdots32$ can be attained in one way and the partial parking outcome $12(n-1)(n-2)\cdots 30$ is attained in $Q_{n-1}$ ways.

		For the partial parking outcome $1(k)(n-1)(n-2)\cdots(k+2)(k+1)0(k-1)(k-2)\cdots32$ for $3\le k\le n-2$, the cars parked after the empty spot must have preferred the spots they ended up in. These cars also have lower indices than every car parked before the empty spot except for car $1$. Thus, no car numbered at least $k$ preferred a spot after the empty spot. If car $1$ also preferred a spot before the empty spot, then the number of ways to attain this partial parking outcome is the number of ways to park car $1$ and cars $k$ through $n-1$ in the order $1(k)(n-1)(n-2)\cdots(k+2)(k+1)$, which is the same as the number of ways to park $n-k+1$ cars in the order $12(n-k+1)(n-k)\cdots43$. Thus, this subcase contributes $Q_{n-k+1}$ to the count. If instead car $1$ preferred a spot after
		the empty spot, there are $k-2$ ways for this to happen, and by the time car $k$ tries to park, car $1$ will have been bumped to spot $1$. Therefore none of the cars from $k$ through $n-1$ prefer spot $1$ and the number of ways for these to park in the order $(k)(n-1)(n-2)\cdots(k+2)(k+1)$ within spots $2$ through $n-k+1$ is the number of ways for $n-k$ cars to park in the order $1(n-k)(n-k-1)\cdots32$, which, by \Cref{cor:difference of adjacent motzkin}, is $M_{n-k+1}-M_{n-k}$. Thus, in total, the partial parking outcome $1(k)(n-1)(n-2)\cdots(k+2)(k+1)0(k-1)(k-2)\cdots32$ for $3\le k\le n-2$ comes from $Q_{n-k+1}+(k-2)(M_{n-k+1}-M_{n-k})$ possible preference lists for the first $n-1$ cars.

		Combining all these cases gives $\left|\Out_{B}^{-1}\left(1(n)(n-1)\cdots32\right)\right|$, which we know, by \Cref{cor:difference of adjacent motzkin}, to be $M_{n+1}-M_{n}$. Hence
		\begin{equation}
			\label{eq:P_n}M_{n+1}-M_{n}=1+2(n-1)+Q_{n-1}+\sum_{k=3}^{n-2}(Q_{n-k+1}+(k-2)(M_{n-k+1}-M_{n-k})).
		\end{equation}
		Reindexing \Cref{eq:P_n} yields
		\begin{equation}
			\label{eq:P_n+1}M_{n+2}-M_{n+1}=1+2n+Q_{n}+\sum_{k=3}^{n-1}(Q_{n-k+2}+(k-2)(M_{n-k+2}-M_{n-k+1})).
		\end{equation}
		Subtracting \Cref{eq:P_n} from \Cref{eq:P_n+1} gives
		\begin{align*}
			M_{n+2}-2M_{n+1}+M_{n} & =2+Q_{n}+\sum_{i=3}^{n-1}(M_{i}-M_{i-1}) \\
			M_{n+2}-2M_{n+1}+M_{n} & =2+Q_{n}+M_{n-1}-M_{2},
		\end{align*}
		and solving for $Q_{n}$ we get $Q_{n}=M_{n+2}-2M_{n+1}+M_{n}-M_{n-1}$, as desired.
	\end{proof}

	\section{Displacement}
	\label{sec:displacement} Given an outcome $\pi=\pi_{1}\pi_{2}\cdots \pi_{n}\in\Sym_{n}$, we recall that $\pi^{-1}=p_{1}p_{2}\cdots p_{n}$ is the permutation of $\{1,2,\dots,n\}$ such that $\pi_{p_j}=j$ for all $j\in [n]$. We recall that the displacement of car $i$ is the absolute difference between the spot car $i$ parks in and the spot it prefers; call this difference $d_{i}$. For a preference list $\alpha$, we define its \textit{displacement multiset} (see \Cref{def:displacementvector}) to be the multiset consisting of the nonzero values $d_{i}$, and we denote this by $\calD_{X}(\alpha)$, where $X\in\{C,M,B\}$ and denotes which parking scheme we are using: $C$ for classical, $M$ for Monopoly, and $B$ for backgammon. We define $N_{X}(\calD_{X})$ to be the number of parking functions under
	the parking scheme $X\in\{C,M,B\}$ with displacement multiset $\calD_{X}$.

	For examples of how displacement multisets differ for classical, Monopoly, and backgammon parking functions, see \Cref{ex:displacement} and \Cref{ex:displacement2}.

	Recall a car is ``lucky'' (otherwise, ``unlucky'') if, at the conclusion of parking, it ends up in the spot that it originally preferred. Let $\lucky(\alpha)$ (resp. $\unlucky(\alpha)$) be the number of lucky (resp. unlucky) cars in the preference list $\alpha$. As we are interested in studying the displacement multiset of Monopoly and backgammon parking functions, we focus on unlucky cars because they contribute to the displacement multiset. While Gessel and Seo \cite{GesselandSeo} established the generating function for the number of classical parking functions with a fixed number of lucky cars, we state their result in terms of unlucky cars:
	\begin{equation}
		\label{eq:GesselSeo}U(q)=\sum_{\alpha\in \PF_n}q^{\unlucky(\alpha)}=\prod_{i=1}^{n-1}(iq+(n-i+1)).
	\end{equation}

	\subsection{Displacement for Monopoly Parking Functions}
	\label{sub:dispMon} Under the Monopoly parking scheme to park $n$ cars, we let $\unlucky_{M}(\alpha)$ denote the number of unlucky cars in the preference list $\alpha\in \MPF_{n}$.

	\begin{theorem}
		\label{thm:repeatedgf} \label{thm:genfunc} The generating function $U_{M}(q)$ for the number of unlucky cars under the Monopoly parking scheme is given by
		\[
			U_{M}(q)=\sum_{\alpha\in \MPF_n}q^{\unlucky_M(\alpha)}=\prod_{i=1}^{n}\left((i-1)q+(n-i+1)\right).
		\]
	\end{theorem}
	\begin{proof}
		Note that car $1$ has $n$ possible preferences, none of which make it unlucky. Parking the cars one at a time, each subsequent car $i$, regardless of how its predecessors have parked, has $i-1$ possible preferences that would make it unlucky (that is selecting a preference among the spots occupied by the previous $i-1$ cars) and thus increase the number of unlucky cars by $1$ and has $n-i+1$ preferences that would not increase the number of unlucky cars. Thus, the generating function is $\Pi_{i=1}^{n}((i-1)q+(n-i+1))$.
	\end{proof}

	Note that the argument in the proof of \Cref{thm:genfunc} would apply to any parking scheme where there are no bumps and every preference list is a parking function. In particular, up to reindexing, this is the same as the analogous generating function given for parking functions of Type B in~\cite{CEFY24}. The same expression will also appear in a different context when we consider the backgammon parking scheme (see Theorem~\ref{thm:backgf}). The coefficients of $q^{i}$ and $q^{n-1-i}$ are always the same and the table of coefficients is given in OEIS \seqnum{A071208}.

	\Monopolyreluctant*
	\begin{proof}
		A displacement multiset of cardinality $1$ means that there is exactly one unlucky car. Specifically, when the displacement multiset is $\{x\}$, there is one doubly preferred spot, which differs by $x$ from the spot that no car preferred. The higher-indexed car with the doubled preference ends up in the undesired spot.

		We begin by counting the number of ways a car that prefers spot $i$ ends up in spot $m$ (and all other cars get what they want). Suppose the unlucky car is car $k$. Its predecessors must use all the spots before spot $m$ but not use spot $m$, which can happen in $\binom{n-m}{k-m}$ ways. If $i<m$, then all these ways have the spot car $k$ wanted already occupied when car $k$ tries to park. If $i>m$, we need to specify that spot $i$ is also used by some predecessor of car $k$, so the predecessors of car $k$ can park in $\binom{n-m-1}{k-m-1}$ ways. There are then $(k-1)!$ ways to permute the preferences of the predecessors and $(n-k)!$ ways to permute the preferences of the successors of car $k$ among the spots still available after car $k$ parks. Now the pairs of $(i,m)$ that give a
		difference of $x$ are $(1,x+1),(2,x+2),\dots(n-x,n)$ and $(x+1, 1), (x+2,2),\dots, (n,n-x)$.

		We will prove the following two claims: \setcounter{claim}{0}
		\begin{claim}
			\label{claimalpha} For $1\le j\le n-x$, the number of ways to get $(j,x+j)$ is $n!/(x+j)$.
		\end{claim}
		\begin{claim}
			\label{claim:two} For $1\le j\le n-x$, the number of ways to get $(x+j,j)$ is $n!/(j+1)$.
		\end{claim}

		\begin{proof}
			[Proof of Claim 1] The number of ways to get $(j,x+j)$ is
			\begin{align}
				\label{case 1}\sum_{k=1}^{n}(k-1)!(n-k)!\binom{n-(j+x)}{k-(j+x)} & =\sum_{k=j+x}^{n}(k-1)!\frac{(n-(j+x))!}{(k-(j+x))!},
			\end{align}
			where, as usual, $\binom{a}{b}=0$ whenever $b<0$.

			Next, we establish the following identity: $\sum_{k=m}^{n}\frac{(k-1)!}{(k-m)!}=\frac{n!}{m(n-m)!}$. We fix $m$ and proceed by induction on $n$. For the base case of $n=m$, we have $(m-1)!=\frac{m!}{m}$. Then, assuming the desired statement holds for $n=n'$, we have
			\begin{align*}
				\sum_{k=m}^{n'+1}\frac{(k-1)!}{(k-m)!} & =\frac{n'!}{(n'+1-m)!}+\sum_{k=m}^{n'}\frac{(k-1)!}{(k-m)!} \\
				                                       & =\frac{n'!}{(n'+1-m)!}+\frac{n'!}{m(n'-m)!}                 \\
				                                       & =\frac{m(n'!)+(n'+1-m)(n'!)}{m(n'+1-m)!}                    \\
				                                       & =\frac{(n'+1)!}{m(n'+1-m)!},
			\end{align*}
			completing the induction. Using this identity in \Cref{case 1} yields that the number of ways to get $(j,x+j)$ is
			\[
				(n-(j+x))!\sum_{k=j+x}^{n}\frac{(k-1)!}{(k-(j+x))!}=(n-(j+x))!\frac{n!}{(j+x)(n-(j+x))!}=\frac{n!}{j+x}.\qedhere
			\]
		\end{proof}

		\begin{proof}
			[Proof of Claim 2]

			Using the same identity, we have that the number of ways to get $(x+j,j)$ is
			\begin{align*}
				\sum_{k=1}^{n}(k-1)!(n-k)!\binom{n-j-1}{k-j-1} & =\sum_{k=j+1}^{n}(k-1)!\frac{(n-j-1)!}{(k-j-1)!} \\
				                                               & =(n-j-1)!\sum_{k=j+1}^{n}\frac{(k-1)!}{(k-j-1)!} \\
				                                               & =(n-j-1)!\frac{n!}{(j+1)(n-(j+1))!}              \\
				                                               & =\frac{n!}{j+1}.
			\end{align*}

			Thus, the number of ways to get $(x+j,j)$ is as claimed.
		\end{proof}

		In total, these cases yield
		\[
			n!\left(\sum_{j=1}^{n-x}\frac{1}{x+j}+\frac{1}{j+1}\right)=n!\left(\sum_{i=x+1}^{n}\frac{1}{i}+\sum_{i=2}^{n-x+1}\frac{1}{i}\right)=n!\left(H_n-H_x+H_{n-x+1}-1\right).
		\]
        \qedhere
	\end{proof}

	\subsection{Displacement for Backgammon Parking Functions}
	\label{sub:dispBack} Under the backgammon parking scheme, the lucky cars are precisely those which are the highest-numbered car with their respective preference. Thus, the number of unlucky cars, $\unlucky_{B}(\alpha)$ in the preference list $\alpha\in \BPF_{n}$ is $n$ minus the number of distinct preferences in $\alpha$.

	We let $\bump(\alpha)$ denote the number of bumps taking place for the preference list $\alpha\in\BPF_{n}$.
	\begin{theorem}
		\label{thm:backgf} The generating function $B(q)$ for the number of bumps under the backgammon parking scheme is given by
		\[
			B(q)=\sum_{\alpha\in\BPF_n}q^{\bump(\alpha)}=\prod_{i=1}^{n}((i-1)q+(n-i+1)).
		\]
	\end{theorem}
	\begin{proof}
		Each car $i$, regardless of how its predecessors have parked, has $i-1$ possible preferences that result in car $i$ performing a bump and $n-i+1$ preferences which do not. Thus car $i$ contributes a factor of $(i-1)q+(n-i+1)$ to the generating function.
	\end{proof}
	Note that this generating function matches the expression in Theorem~\ref{thm:repeatedgf} for the number of unlucky cars when parking under the Monopoly parking scheme. We now give a count for backgammon parking functions with a given displacement multiset with one element, analogous to the Monopoly result presented in \Cref{thm:Monopoly1reluctant}. \backgammondisplacement*

	\begin{proof}
		Exactly $n-1$ cars park in their preferred spot. The remaining car must be in the one spot that was undesired, which means the duplicated preference must differ by $x$ from the label of the undesired spot. Under the backgammon parking scheme, the one car that gets bumped is the lower-numbered car with the duplicated preference and since every other car ends up in its preferred spot, the bumped car automatically ends up in the undesired spot. That is, only the multiset of preferences matters, not the order. There are $2(n-x)$ ways to pick the pair of the duplicate preference and the undesired spot. There are $\binom{n}{2}$ ways to pick which two cars have the same preference and $(n-2)!$ ways to permute the preferences of the other cars. This gives a total of $(n-x)n!$ backgammon
		parking functions with displacement multiset $\{x\}$.
	\end{proof}

	We now consider the number of backgammon parking functions with a displacement multiset with two elements. \backgammondisplacementtwo*

	\begin{proof}
		For a displacement multiset $\{x_{1},x_{2}\}$ with $x_{1}\ge x_{2}$, we can have either three cars with a preference for the same parking spot, while all other cars have a distinct preference, or we can have two pairs of cars with the same preference (but distinct from each other). The first case we refer to as having ``a triple-preference'', and the second case we refer to as having ``two double-preferences''. Note that in either case, the unlucky cars end up in the two undesired spots, so their positions (up to swapping which unlucky car ends up in which undesired spot) are determined by the multiset of preferences. We now consider each case separately and start by letting $T(x_{1},x_{2})$ denote the number of preference lists arising from the triple-preference case.

		To find the value of $T(x_{1},x_{2})$, we first count preference lists containing a triple-preference for spot $p$ which have no other repeated preferences and satisfy either of the following conditions:
		\begin{itemize}
			\item neither spot $p+x_{1}$ nor $p-x_{2}$ is preferred by any car,

			\item neither spot $p-x_{1}$ nor $p+x_{2}$ is preferred by any car.
		\end{itemize}

		If $x_{1}+x_{2}\geq n$, then no such triple-preferences can occur as the undesired spots $p-x_{1},p+x_{2}$ or $p-x_{2},p+x_{1}$ would differ by at least $n$. Hence, in this case, it must be that $x_{1}+x_{2}<n$. For $x_{1}\ne x_{2}$, there are $2(n-x_{1}-x_{2})$ such triple-preferences as $p$ can range from $1+x_{1}$ to $n-x_{2}$ or from $1+x_{2}$ to $n-x_{1}$ depending on the subcase. For $x_{1}=x_{2}$, there are only $n-x_{1}-x_{2}=n-2x_{1}$ such triple-preferences as the pairs $p-x_{1},p+x_{2}$ and $p-x_{2},p+x_{1}$ now coincide and $p$ can range from $1+x_{1}$ to $n-x_{1}$. If $x_{1}> x_{2}$, we also have to consider when the undesired spots are both higher or both lower than the triple-preference $p$, which yields another $2(n-x_{1})$ possible choices, as $p$ can range from $1$
		to $n-x_{1}$ when both undesired spots are greater than $p$ or from $x_{1}+1$ to $n$ when both undesired spots are less than $p$. Then there are $\binom{n}{3}$ ways to pick which cars have the exact same preference, times $(n-3)!$ ways to decide which remaining cars get which remaining preferences. Thus, the total number of preference lists with a triple-preference and displacement multiset $\{x_{1},x_{2}\}$ is

		\[
			T(x_{1},x_{2})=
			\begin{cases}
				\displaystyle \frac{n!}{6}\bigl(4n-4x_{1}-2x_{2}\bigr), & \text{ for }x_{1}>x_{2},\ n\ge x_{1}+x_{2}, \\[1.2em] \displaystyle \frac{n!}{6}\bigl(2n-2x_{1}\bigr), & \text{ for } x_{1}>x_{2},\ n<x_{1}+x_{2}, \\[1.2em] \displaystyle \frac{n!}{6}\bigl(n-2x_{1}\bigr) , & \text{ for } x_{1}=x_{2},\ n\ge 2x_{1}, \\[1.2em] \displaystyle 0, & \text{ for }x_{1}=x_{2},\ n<2x_{1}.
			\end{cases}
		\]

		In the case of two double-preferences, the lower-numbered car in each double ends up in an undesired spot. If the double-preferences are spots $a$ and $b$ (with $a\neq b$), the undesired spots are $c$ and $d$ with $|a-c|=x_{1}, |b-d|=x_{2}$. We have to consider the case where the unlucky car that preferred spot $a$ parks in spot $d$ and where the unlucky car that preferred spot $b$ parks in spot $c$. If this gives the same displacement multiset, it means that either $|a-d|=x_{1}, |b-c|=x_{2}$ or that $|a-d|=x_{2}, |b-c|=x_{1}$. In the former case, we get that $a$ is equidistant from $d$ and $c$ so it is exactly in between, but $b$ is also equidistant from $c$ and $d$, which gives a contradiction. In the latter case, $d$ is equidistant from $a$ and $b$, and $c$ is also equidistant
		from $a$ and $b$, again giving a contradiction. Thus, we only get the displacement multiset $\{x_{1},x_{2}\}$ if the unlucky car which preferred spot $a$ parks in spot $c$ and unlucky car which preferred spot $b$ parks in spot $d$.

		Next, we determine the number of ways to select the double-preferred spots $a,b$ and the undesired spots $c$ and $d$. The pair $a,c$ differs by $x_{1}$ and the pair $b,d$ differs by $x_{2}$. There are $n-x_{1}$ ways to pick a pair that differs by $x_{1}$ and $n-x_{2}$ ways to pick a pair that differs by $x_{2}$, but we also have to make sure not to count the cases where these pairs overlap.

		When $x_{1}\ne x_{2}$, the overlap occurs when some spot is $x_{1}$ away from the second and $x_{2}$ away from the third, which is precisely what we tabulated previously when counting triples.

		This leaves $(n-x_{1})(n-x_{2})-(4n-4x_{1}-2x_{2})=n^{2}-(x_{1}+x_{2}+4)n+(x_{1}x_{2}+4x_{1}+2x_{2})$ possibilities when $x_{1}+x_{2}\le n$ and $(n-x_{1})(n-x_{2})-2(n-x_{1})=n^{2}-(x_{1}+x_{2}+2)n+(x_{1}x_{2}+2x_{1})$ possibilities when $x_{1}+x_{2}>n$.

		When $x_{1}=x_{2}$, we have $\frac{1}{2}(n-x_{1})(n-x_{1}-1)$ ways to pick two distinct pairs with each differing by $x_{1}$, but for $n\ge 2x_{1}$, we must subtract the $n-2x_{1}$ ways where they overlap. This leaves $\frac{1}{2}(n^{2}-(2x_{1}+3)n+x_{1}^{2}+5x_{1})$ options for $n\ge 2x_{1}$ and $\frac{1}{2}(n-x_{1})(n-x_{1}-1)$ for $n<2x_{1}$.

		Since the pair $a,c$ differs by $x_{1}$ and the pair $b,d$ differs by $x_{2}$, the freedom to interchange the values of $a$ and $c$ contributes a factor of 2 to the count. Similarly, we get a factor of 2 for the pair $b, d.$

		Hence, the number of ways to select the double-preferred spots $a,b$ and the undesired spots $c$ and $d$ is given by:
		\[
			\begin{cases}
				\displaystyle 4\bigl(n^{2}-(x_{1}+x_{2}+4)n +(x_{1}x_{2}+4x_{1}+2x_{2})\bigr), & \text{ for }x_{1}>x_{2},\ n\ge x_{1}+x_{2}, \\

				\displaystyle 4\bigl(n^{2}-(x_{1}+x_{2}+2)n +(x_{1}x_{2}+2x_{1})\bigr), & \text{ for }x_{1}>x_{2},\ n<x_{1}+x_{2}, \\[1em]

				\displaystyle 2\bigl(n^{2}-(2x_{1}+3)n+x_{1}^{2}+5x_{1}\bigr), & \text{ for }x_{1}=x_{2},\ n\ge 2x_{1}, \\[1.2em]

				\displaystyle 2\bigl(n^{2}-(2x_{1}+1)n+(x_{1}^{2}+x_{1})\bigr), & \text{ for }x_{1}=x_{2},\ n<2x_{1}.
			\end{cases}
		\]

		In addition, there are $\binom{n}{2}\binom{n-2}{2}$ ways to pick which cars prefer spots $a,b$ and $(n-4)!$ ways to assign the remaining cars to the remaining preferred spots. For bookkeeping purposes, we denote this count by:
		\[
			K(x_{1},x_{2}):=
			\begin{cases}
				\displaystyle n!\bigl(n^{2}-(x_{1}+x_{2}+4)n +(x_{1}x_{2}+4x_{1}+2x_{2})\bigr), & \text{ for }x_{1}>x_{2},\ n\ge x_{1}+x_{2}, \\

				\displaystyle n!\bigl(n^{2}-(x_{1}+x_{2}+2)n +(x_{1}x_{2}+2x_{1})\bigr), & \text{ for }x_{1}>x_{2},\ n<x_{1}+x_{2}, \\[1em]

				\displaystyle \frac{n!}{2}\bigl(n^{2}-(2x_{1}+3)n+x_{1}^{2}+5x_{1}\bigr), & \text{ for }x_{1}=x_{2},\ n\ge 2x_{1}, \\[1.2em]

				\displaystyle \frac{n!}{2}\bigl(n^{2}-(2x_{1}+1)n+(x_{1}^{2}+x_{1})\bigr), & \text{ for }x_{1}=x_{2},\ n<2x_{1}.
			\end{cases}
		\]
		We claim that exactly half of these $K(x_{1},x_{2})$ possible preference lists actually give displacement multiset $\{x_{1},x_{2}\}$. Assuming this claim holds, the total number of preference lists that produce the displacement multiset $\mathcal{D}_{B}=\{x_{1},x_{2}\}$ is then $N_{B}(\{x_{1},x_{2}\})=T(x_{1},x_{2})+\frac{K(x_1,x_2)}{2},$ which simplifies to

		\begin{align}
			\label{eq:T +K/2}\raisetag{10ex}{\small N_B(\{x_1,x_2\})= \begin{cases}\displaystyle \frac{n!}{6}\bigl( 3n^{2} -(3x_{1}+3x_{2}+8)n +(3x_{1}x_{2}+8x_{1}+4x_{2}) \bigr), & \text{ for }x_{1}>x_{2},\ n\ge x_{1}+x_{2}, \\[1.2em] \displaystyle \frac{n!}{6}\bigl( 3n^{2} -(3x_{1}+3x_{2}+4)n +3x_{1}x_{2}+4x_{1} \bigr), & \text{ for } x_{1}>x_{2},\ n<x_{1}+x_{2}, \\[1.2em] \displaystyle \frac{n!}{12}\bigl( 3n^{2} -(6x_{1}+7)n +(3x_{1}^{2}+11x_{1}) \bigr), & \text{ for } x_{1}=x_{2},\ n\ge 2x_{1}, \\[1.2em] \displaystyle \frac{n!}{4}\bigl( n^{2} -(2x_{1}+1)n +(x_{1}^{2}+x_{1}) \bigr), & \text{ for }x_{1}=x_{2},\ n<2x_{1}.\end{cases} }
		\end{align}
		Notice that Equation \eqref{eq:T +K/2} recovers Equations \eqref{eq:x1 neq x2} and \eqref{eq:x1 eq x2}, which completes the proof of the theorem.

		We now prove the claim. If the four cars with preferred spots $a,b$ are numbered $y_{1}<y_{2}<y_{3}<y_{4}$, then in the partial parking outcome after step $y_{4}$, there are two cars not in their preferred spots and the relative order is determined by the preference of car $y_{4}$.

		If car $y_{4}$ prefers spot $b$, then the most recent car to be bumped is the unlucky car preferring spot $b$. Thus, the unlucky car preferring spot $b$ occupies a higher-numbered spot in the partial parking outcome after step $y_{4}$ than the unlucky car preferring spot $a$. Similarly, if $y_{4}$ prefers spot $a$, then it will bump the unlucky car preferring spot $a$ into a higher-numbered spot than the unlucky car preferring spot $b$. Thus, in the partial parking outcome after step $y_{4}$, the two possible relative orders of the unlucky cars preferring spots $a$ and $b$ are equally likely.

		For a fixed set of preferences for the other $n-4$ cars, the relative order of the unlucky cars in the partial parking outcome after step $y_{4}$ determines the final relative order of the unlucky cars, namely which ends up in spot $c$ and which ends up in spot $d$. Thus, outcomes where the unlucky cars which preferred spots $a$ and $b$ park in spots $c$ and $d$, respectively, and outcomes where they park in spots $d$ and $c$, respectively, are equally likely. Thus, the total number of preference lists with displacement multiset $\{x_{1},x_{2}\}$ arising from the two double-preferences case is $\frac{1}{2}K(x_{1},x_{2})$, which completes the proof of the claim. \qedhere
	\end{proof}
	\subsection{A Comparison}
	\label{sub:dispComp} We conclude this section by comparing the number of backgammon and Monopoly parking functions for the same displacement multiset with one element. \displacementcomparison*

	\begin{proof}
		For there to be exactly one unlucky car, there must be two cars with the same preference, while every other car has a distinct preference. For the one unlucky car to have displacement $x$, the doubly preferred spot and the undesired spot must differ by $x$. Under the backgammon parking scheme, any such preference list yields displacement multiset $\{x\}$, but under the Monopoly parking scheme, there is a possibility that the higher-indexed car with the repeated preference blocks a later car from its preferred spot, causing another car to be unlucky. Thus,
		\begin{align}
			\label{inequality in theorem}N_{M}(\{x\}) & \le N_{B}(\{x\}).
		\end{align}
		When $n=2$, the only option is that $x=1$. In this case, there are only two spots and only two cars, so the only such preference lists are $(1,1)$ and $(2,2)$ both of which have displacement multiset $\{1\}$ for both the Monopoly and backgammon parking schemes. Thus, in this case, we achieve equality in \eqref{inequality in theorem}.

		When $n\ge 3$, regardless of the value of $x$, there is some preference list that has displacement multiset $\{x\}$ under the backgammon parking scheme, but not under the Monopoly parking scheme. In particular, when $x=1$, consider any preference list where the first two cars prefer spot $2$ and the remaining cars have distinct preferences in $[n]\setminus\{2,3\}$. When $x>1$, consider any preference list where the first two cars prefer spot $1$ and the remaining cars have distinct preferences in $[n]\setminus\{1,x+1\}$. Thus, in this case, the inequality in \eqref{inequality in theorem} is strict. \qedhere
	\end{proof}

	\section{Future directions}
	\label{sec:future} We conclude this work by providing some questions for further study.

	In Theorem \ref{thm:at most Catalan}, we established that the maximum value of $|\Out_{B}^{-1}(\pi)|$ is the Catalan number $C_{n}$, while in Theorem \ref{thm:backgammon lower bound}, we established that the minimum value of $|\Out_{B}^{-1}(\pi)|$ is at least $2^{n-1}$. In Table \ref{tab:backmin}, we note the minimum values of $|\Out_{B}^{-1}(\pi)|$ computed for up to ten cars, along with the outcomes $\pi$ obtaining these minima. The sequence does not appear on the OEIS, hence, we ask the following:

	\begin{question}
		For $\pi\in\mathcal{\Sym}_{n}$, what is the minimum value of $|\Out_{B}^{-1}(\pi)|$? Moreover, which outcomes achieve the minimum value?
	\end{question}

	\begin{center}
		\begin{table}[h]
			\centering
			\begin{tabular}{c|c}
				Minimum $|\Out_{B}^{-1}(\pi)|$ & Outcomes achieving minimum   \\
				\hline
				1                              & 1                            \\
				2                              & 12, 21                       \\
				4                              & 213, 231, 321                \\
				8                              & 3214, 3241                   \\
				17                             & 43215, 43251, 43521          \\
				36                             & 435216, 435261               \\
				75                             & 5463217, 5463271             \\
				160                            & 65473218, 65473281           \\
				351                            & 654873219, 654873291         \\
				726                            & 765984321(10), 76598432(10)1
			\end{tabular}
			\caption{Minimum value of $|\Out_{B}^{-1}(\pi)|$ for up to ten cars and each outcome $\pi$ obtaining the minimum.}
			\label{tab:backmin}
		\end{table}
	\end{center}

	While it remains an open problem to determine both the minimum cardinality of $|\Out_{B}^{-1}(\pi)|$, as well as exactly which outcomes $\pi$ achieve that minimum, the data in \Cref{tab:backmin} suggests there may be a pattern characterizing the outcomes that achieve the minimum value. We leave this exploration to future work.

	In \Cref{thm: everything is Monopoly} we determined exactly the number of preference lists yielding each possible outcome under the Monopoly parking scheme. However, we only determined exact values under the backgammon parking scheme for certain outcomes. For instance, in Theorem \ref{thm:motzkin} and its corollaries, we determined $|\Out_{B}^{-1}((n)(n-1)\cdots21)|, |\Out_{B}^{-1}(1(n)(n-1)\cdots2)|$, and $|\Out_{B}^{-1}(12(n)(n-1)\cdots 3)|$, all of which can be expressed as a sum involving Motzkin numbers. Hence, we ask:

	\begin{question}
		For fixed $k\in\mathbb{N}$, can $|\Out_{B}^{-1}(12\cdots(k-1)(k)(n)(n-1)\cdots(k+2)(k+1))|$ for $n>k$ be expressed as a sum involving Motzkin numbers whose indices depend on $n$?
	\end{question}

	Next, we provide some questions related to displacement. Recall that we studied the number of Monopoly and backgammon parking functions with a prescribed displacement multiset of small cardinality. Hence, it is natural to ask the following:

	\begin{question}
		Determine counts or bounds for the numbers $N_{M}(Y)$ and $N_{B}(Y)$ for multisets $Y$.
	\end{question}

	In \Cref{thm:displacementcomparison}, we showed that $|N_{M}(Y)|\le |N_{B}(Y)|$ when $|Y|=1$. It would also be interesting to classify for which other $Y$ this inequality holds.

	Given a fixed preference list $\alpha$, the partial parking outcomes $\Out_{M,t}(\alpha)$ and $\Out_{B,t}(\alpha)$ under the Monopoly and backgammon parking schemes have the same spots filled after every step $t$. However, they often lead to distinct displacement multisets. Along these lines, we ask:

	\begin{question}
		For $n$ cars, can we characterize and enumerate the preference lists that have the same displacement multiset under the backgammon and Monopoly parking schemes?
	\end{question}

	One can always consider a street with more spots than cars; hence we ask:
	\begin{question}
		How do our results on backgammon and Monopoly parking functions change if we are allowed to consider more spots than cars?
	\end{question}

	\section*{Acknowledgments}
	P.\ E.\ Harris was supported in part by a Simons Travel Support for Mathematicians award \#00007621. A.\ Clifton was supported in part by the Institute for Basic Science (IBS-R029-C1) and by the grant 23-06815M of the Grant Agency of the Czech Republic. The contents of this work were solely generated by humans; no AI was used.

	\bibliographystyle{plain}
	\bibliography{Bibliography.bib}
\end{document}